\documentclass[11pt]{amsart}

\def\squarebox#1{\hbox to #1{\hfill\vbox to #1{\vfill}}}

\usepackage{amssymb}
\usepackage{amsmath}
\usepackage{epsfig}

\newcommand{\bel}{\begin{equation} \label}
\newcommand{\ee}{\end{equation}}

\newcommand{\R}{{\mathbb R}}
\newcommand{\RR}{{\mathbb R}}

\newcommand{\C}{{\mathbb C}}
\newcommand{\N}{{\mathbb N}}
\newcommand{\Tr}{\operatorname{Tr}}

\newcommand{\rank}{\operatorname{rank}}
\newcommand{\mult}{\operatorname{mult}}

\newcommand{\supp}{\operatorname{supp}}

\newcommand{\im}{\operatorname{Im}}

\newcommand{\loc}{\operatorname{loc}}

\newcommand{\tr}{{\operatorname{tr}}}

\renewcommand{\Re}{\mathop{\rm Re}\nolimits}
\renewcommand{\Im}{\mathop{\rm Im}\nolimits}

\newcommand{\fract}[2]{\genfrac{}{}{0pt}{}{\scriptstyle #1}{\scriptstyle #2}}
\theoremstyle{plain}

\newtheorem{thm}{Theorem}
   
\newtheorem{cor}{Corollary}
   
\newtheorem{lem}{Lemma}
   
\newtheorem{prop}{Proposition}
   
\newtheorem{defi}{Definition}

\theoremstyle{definition}
\newtheorem*{rem}{Remark}
\newtheorem*{rems}{Remarks}

\numberwithin{equation}{section}

\title[Resonances and SSF near the  Landau levels]
{Resonances and Spectral Shift Function\\ near the Landau levels}
\author[J.F. Bony, V. Bruneau]{Jean-Fran\c{c}ois Bony, Vincent Bruneau }
\address {Institut de Math\'ematiques de Bordeaux, FR CNRS 2254, MAB UMR CNRS 5466
Universit\'e Bordeaux I,  351, Cours de la Lib\'eration, 33405  Talence, France}\email{bony@math.u-bordeaux1.fr, vbruneau@math.u-bordeaux1.fr}
\author[G. Raikov]{Georgi Raikov}
\address{Departemento de Mathem\'aticas, Facultad de Ciencas, Universidad de Chile, Las Palmeras 3425, Santiago, Chile}
\email{graykov@uchile.cl}

\def\phi {\varphi}

\newcommand{\vp}{v_\perp}
\newcommand{\xp}{X_\perp}
\newcommand{\mper}{m_\perp}

\newcommand{\res}{{\rm Res}(H)}

\begin{document}

\begin{abstract}
We consider the 3D Schr\"odinger operator $H = H_0  + V$ where $H_0 =
(-i\nabla - A)^2 -b$, $A$ is a magnetic potential generating a constant magnetic
field of strength $b>0$, and $V$ is a short-range electric potential which
decays superexponentially with respect to the variable along the magnetic
field. We show that the resolvent of $H$ admits a meromorphic extension from
the upper half plane to an appropriate Riemann surface ${\mathcal M}$, and
define the resonances of $H$ as the poles of this meromorphic extension. We
study their distribution near any fixed Landau level $2bq$, $q \in {\mathbb
  N}$.
First, we obtain a sharp upper bound of the number of resonances
in a vicinity of $2bq$. Moreover, under appropriate hypotheses, we
establish corresponding lower bounds which imply the existence of
an infinite number of resonances, or the absence of resonances in
certain sectors adjoining $2bq$. Finally, we deduce a
representation of the derivative of the spectral shift function
(SSF) for the operator pair $(H,H_0)$ as a sum of a harmonic
measure related to the resonances, and the imaginary part of a
holomorphic function. This representation justifies the
Breit-Wigner approximation, implies a trace formula, and provides
information on the singularities of the SSF at the Landau levels.
\end{abstract}

\maketitle

\bigskip

\bigskip

\section{Introduction}

Let
$$
H_0: =
(D_{x_1} + \frac{b}{2} x_2)^2 + (D_{x_2} - \frac{b}{2} x_1)^2 -b +
D_{x_3}^2 ,\quad D_\nu := -i \frac{\partial}{\partial \nu},
$$
be the Schr\"odinger operator with homogeneous magnetic field of strength
$b>0$,
pointing at the $x_3$-direction.  Initially, the self-adjoint operator $H_0$ is
defined on $C^{\infty}_0(\R^3)$, and then is closed in $L^2(\R^3)$.
 This operator can be  written in  $L^2(\R^3) = L^2(\R^2) \otimes L^2(\R)$
  as
$$H_0=H_{0, \perp} \otimes I + I \otimes D_{x_3}^2, $$
with $H_{0, \perp}=
(D_{x_1} + \frac{b}{2} x_2)^2 + (D_{x_2} - \frac{b}{2} x_1)^2-b $.

It is
well known that the spectrum of the operator $H_{0, \perp}$  consists of the Landau levels
 $2qb$,  $q \in \N : = \{0,1\ldots\}$, and the multiplicity of each
eigenvalue $2bq$ is infinite (see e.g. \cite{ahs}). Consequently,
the spectrum of $H_0$ is absolutely continuous, equals $\lbrack
0,+\infty\lbrack$, and has an infinite set of thresholds $2qb$,
$q\geq 0$.

For ${\bf x} = (x_1,x_2,x_3) \in \RR^3$ we denote by $\xp =
(x_1,x_2)$ the variables on the plane perpendicular to the
magnetic field.  We assume that the electric potential  $V: \R^3
\to \R$ is Lebesgue measurable, and satisfies the estimates
\begin{equation} \label{eq1.1}
 V({\bf x})={\mathcal O}( \langle \xp \rangle^{-\mper}\;\langle  x_3 \rangle^{-m_3}),
\quad {\bf x} \in \RR^3,
 \end{equation}
with  $\mper >2$,  $m_3 >1$, and $\langle  x \rangle : =
(1 + | {x}|^2)^{1/2}$, $ x \in \RR^d$, $d \geq 1$.



 On the domain of $H_0$ we introduce the operator
$H :=H_0 + V$.
  Since $V$ is a relatively compact perturbation of $H_0$, it follows from the
 Weyl criterion that the essential spectra of
$H$ and $H_0$ are the same.  Moreover, since $V$ is a relatively trace-class
 perturbation, the Kato-Rosenblum theorem implies that the absolutely
 continuous spectrum coincides with $[0,+\infty[$.

 It is well  known that $H$ can  have infinite  negative
discrete spectrum  and, for some special $V$, it can have
 infinitely  many embedded eigenvalues below each Landau
level (see \cite{ahs}, \cite{RCubo} or \cite{RGiens}).   On the other
hand, it is shown in \cite{FR} that in the case of sign-definite
$V$, the spectral shift function (SSF)  for the operator pair $(H,H_0)$ has a singularity at
each Landau level.  Therefore, it is natural to  expect
that  there could be accumulation of the resonances of the operator $H$ near the
Landau levels.  For a Coulomb potential, some numerical
results confirm this conjecture \cite{Resnum}.   The
goal of this paper is to study the resonances near the Landau levels,
and to establish the link between these resonances and the
spectral shift function by the so-called Breit-Wigner
approximation.  Such a representation of the derivative of
the spectral shift function related to the resonances, implies
trace formulas which have given recently a substantial impetus to
the research concerning the upper and lower bounds of the number
of resonances in different situations (see \cite{SjBook}, \cite{S1}, \cite{S3},  \cite{S2}, \cite{PZ3}, \cite{BSj}, \cite{BP},  \cite{DZ}).

 We consider potentials $V$ which decay super-exponentially with respect to
 $x_3$ (or  are compactly supported with respect to $x_3$). Hence, we
do not use
 dilation methods in order to
define the resonances near the Landau levels.  For a
definition using
 complex dilation, we refer  the
 reader to  \cite{Wa}, \cite{FrWa} where precise asymptotics as $b \to \infty$
of the resonances near the real axis is given.
In our work, $b$ is fixed, and we study the number of resonances in a domain
 $2bq + r \Omega$ as $r$ tends to $0$. Then we justify the Breit-Wigner
 approximation for the spectral shift function near  the Landau
 levels.

The paper is organized as follows. In the next  section,  we
define  the resonances as the poles of the resolvent, the
first step being to introduce   a Riemann surface to which the
resolvent is extended. Note that the resonances defined as poles
of the resolvent, are also zeros of  a generalized Krein
perturbation determinant  with the same multiplicity. In  Section 3, we  obtain an
upper bound of the number of resonances in a domain $2bq + r
\Omega$ as $r$ tends to $0$ (see Theorem \ref{thmUB}).  In
 Section 4,  we obtain more  information on the
localization of the resonances for the  case of perturbations of
definite sign.  In particular, we show that there is an infinite
number of resonances near any arbitrary fixed  Landau level  for
small $V$ of sufficiently rapid decay (see Theorem \ref{smallV}),
and that there are no embedded eigenvalues for small positive $V$
(see  Proposition \ref{cor2}).  At last, in Section 5, we
represent  the derivative of the spectral shift function
near  the Landau levels as a sum of  a harmonic measure
related to the resonances and the imaginary part of a holomorphic
function (see Theorem \ref{thmBW}).  Such a representation justifies
 the Breit-Wigner approximation, implies  a trace
formula, and  for a special class of $V$ sufficiently slowly decaying with respect to
the
variables perpendicular to the magnetic field, allows us to
estimate the remainder in the asymptotic relations obtained in
\cite{FR}.









\section{Resonances}
In this section we define the resonances of $H=H_0 + V$ for $V$
decaying super-exponentially with respect to $x_3$,  i.e.
\begin{equation}  \label{supexp}
V({\bf x}) = {\mathcal O} ( \langle \xp \rangle^{-\mper} \exp(-N |x_3|)),
\end{equation}
for $\mper \geq 0$ and any $N>0$. As in \cite{SaZw}, the
resonances will be defined as the poles of the meromorphic
continuation of the resolvent in some weighted $L^{2}$ spaces.
Since $V$ is not compactly supported with respect to $x_3$, the
cut-off resolvent cannot be used here.

First, we have to prove  the existence of a holomorphic
extension for the unperturbed operator.  Let ${\mathbb C}_+
: = \{\lambda \in {\mathbb C}; \ {\rm Im}\, \lambda > 0\}$ be the open
upper half plane.  For $\lambda \in {\mathbb C}_+$ we have
\begin{equation}\label{eq2.1}
(H_0- \lambda )^{-1}= \sum_{q \in \N} p_q \otimes (D_{x_3}^2 + 2bq
- \lambda )^{-1},
\end{equation}
where $p_q$ is the orthogonal projection onto  ${\mathcal H}_q:=
{\hbox{ker}}(H_{0, \perp}-2bq)$.

Let us recall that for $\Im k>0$, the integral kernel of the
operator $(D^{2}_{x_{3}} - k^{2})^{-1}$ is given by
\begin{equation}  \label{mml}
{\mathcal R}(x_{3} , x_{3}' )= \frac{i e^{i k |x_{3} -x_{3}' |}}{2 k }.
\end{equation}
Then, for $N>0$, the  operator-valued function
\begin{equation}  \label{tchi}
t_{N} (k^2) := e^{-N \langle x_{3} \rangle} (D_{x_3}^2 -k^2)^{-1} e^{-N
  \langle x_{3} ' \rangle} \in {\mathcal L}(L^{2} (\R_{x_{3}}), H^{2}
  (\R_{x_{3}})) ,
\end{equation}
can be extended holomorphically from ${\mathbb C}_+$ to  $\{ k \in \C^{*} ; \ {\rm Im}
\, k > -N \}$.  Hence, for any $N>0$ and $q \in \N$,
\begin{align*}
z \mapsto (D_{x_3}^2 + 2bq -z)^{-1}  \in {\mathcal L}(e^{-N \langle
  x_{3} \rangle} L^2
(\R_{x_3}) , e^{N \langle x_{3} \rangle} H^2 (\R_{x_3})) ,
\end{align*}
has a holomorphic extension from $\C \setminus [2bq, + \infty[$ to the
$2$-sheeted covering $\pi_q: k\in \C^* \mapsto  k^2+ 2bq \in \C \setminus \{ 2b q\}$
(with $\im k >-N$). However,
this covering depends on $q$, and therefore it is not suitable for the
extension of
\eqref{eq2.1}.

A natural domain of analytic extension of \eqref{eq2.1} is the  universal covering of $\C\setminus 2b
\N$:
$$\overline{ \pi} : \overline{\C\setminus 2b \N} \rightarrow \C\setminus 2b
\N,$$
but it  does not give a maximal analytic continuation. Indeed for some $z
\in \C\setminus 2b \N$, there are some points $z_1$, $z_2 \in\overline{
  \pi}^{-1}(z)$, $z_1\neq z_2$ such that the germs of $(H_0-z)^{-1}$ at $z_1$
and at $z_2$ are the same.  Next, we introduce the equivalence relation
${\mathcal R}$ concerning such pairs of points $(z_1,z_2)$. Let $\pi_1(\C\setminus 2b
\N)$ be the fundamental group of $\C\setminus 2b \N$, and $G$ be its subgroup
generated by $\{a_1^2, \; a_2 a_1 a_2^{-1}a_1^{-1}; \; \hbox{with
} a_1, a_2 \in \pi_1(\C\setminus 2b \N)\}$.  We will write $z_1
{\mathcal R} z_2$ if and only if $ \overline{ \pi}z_1=\overline{
\pi}z_2$, and for any path $\gamma$  connecting $z_1$ with $z_2$,
the class of the closed path $\pi \gamma$ in $\pi_1(\C\setminus 2b
\N)$ is an element of $G$
(that is $\pi(\gamma)$  goes an even number of  times round each
Landau level).

Then we define the domain $ {\mathcal M}$ of the analytic
extension of $(H_0-z)^{-1}$ as the quotient of $
\overline{\C\setminus 2b \N}$ by the relation ${\mathcal R}$. This
domain can also be identified with the following covering of
$\C\setminus 2b \N$ (see for instance  Proposition 13.23 of
\cite{Fulton}):

\begin{defi}
Let $\pi_1(\C\setminus 2b \N)$  be the fundamental group of
$\C\setminus 2b \N$.  Let $G$ be the subgroup of
$\pi_1(\C\setminus 2b \N)$ generated by $\{a_1^2, \;a_2 a_1 a_2^{-1}a_1^{-1}; \; \hbox{with }
a_1, a_2 \in \pi_1(\C\setminus 2b \N)\}$.
We define $\pi_G: \; {\mathcal M} \rightarrow \C\setminus 2b \N$ as
the connected infinite-sheeted covering such that
$\pi_1({\mathcal M})=G$.
\end{defi}

From now on, we fix a base point in ${\mathcal M}$, and define
the physical plane ${\mathcal F}$ as the connected
component of $\pi_{G}^{-1} (\C \setminus [0, + \infty[)$ containing this
base point. By definition, the functions ${\mathcal M} \ni z
 \mapsto \sqrt{z-2b q}$ have a positive imaginary part on
${\mathcal F}$. Let ${\mathcal F}_{+} = {\mathcal F} \cap \pi_{G}^{-1}
(\C_{+})$ be the upper half-plane.
In what follows, we identify
${\mathcal F}$ (resp. ${\mathcal F}_{+}$ and $\partial
{\mathcal  F}_{+}$) with $\C \setminus [0, + \infty[$ (resp. $\C_{+}$ and $\R \setminus 2
b \N$), and  denote by $z$ the generic point on ${\mathcal M}$.

For $\lambda_0 \in {\mathbb C}$ and $\varepsilon
  > 0$ put $D(\lambda_0, \varepsilon) : = \{\lambda \in {\mathbb C}, |\lambda
  - \lambda_0| < \varepsilon\}$  and $D(\lambda_0, \varepsilon)^{*}: = \{\lambda \in {\mathbb C}, 0<|\lambda
  - \lambda_0| < \varepsilon\}$. 

\begin{defi}  \label{a36}
We denote by $D_{q}^{*} \subset {\mathcal M}$, the connected component of
 $\pi_{G}^{-1} ( D (2b q, 2b)^{*})$ that intersects
${\mathcal F}_{+}$.

Since   $\pi_{G}:D_{q}^{*} \to D (2b q, 2b)^{*}$ is a
2-sheeted covering of $D (2b q, 2b)^{*}$, there
exists an analytic bijection
\begin{equation}  \label{a35}
z_{q} : k \in D (0 , \sqrt{2b})^{*}\to z_{q} (k) \in
D_{q}^{*} ,
\end{equation}
such that $\pi_{G} (z_{q} (k)) = 2 b q + k^{2}$ and $z_{q}^{-1}
(D_{q}^{*} \cap {\mathcal F}_{+} )$ is the first quadrant of  $D (0 ,
\sqrt{2b})^{*}$.

For $N>0$, we denote by
${\mathcal M}_{N}$ the set of points $m \in {\mathcal M}$ such
that for each $q \in \N$, we have  $\Im \sqrt{z-2b q} > -N$. Of course, we
have $\cup_{N>0} {\mathcal M}_{N} = {\mathcal M}$.
\end{defi}

%
%
%
%

Figure 1 summarizes the setting near the Landau level $2b q$.
For the free operator, we have the following proposition.

\begin{figure} \label{a100}
\begin{center}
\begin{picture}(0,0)%
\includegraphics{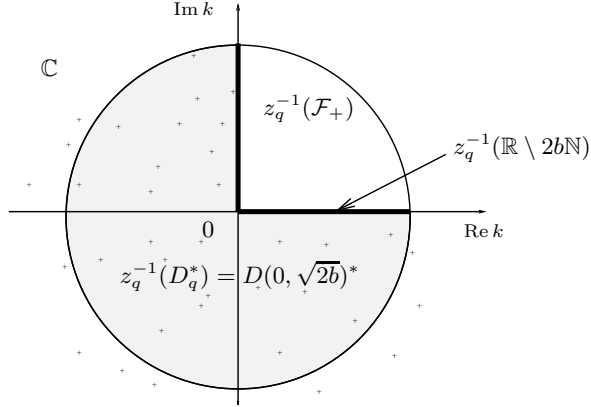}%
\end{picture}%
\setlength{\unitlength}{1184sp}%
\begingroup\makeatletter\ifx\SetFigFont\undefined%
\gdef\SetFigFont#1#2#3#4#5{%
  \reset@font\fontsize{#1}{#2pt}%
  \fontfamily{#3}\fontseries{#4}\fontshape{#5}%
  \selectfont}%
\fi\endgroup%
\begin{picture}(10041,8496)(1168,-8044)
\put(5251,-4486){\makebox(0,0)[lb]{\smash{{\SetFigFont{9}{10.8}{\rmdefault}{\mddefault}{\updefault}$0$}}}}
\put(3526,-5461){\makebox(0,0)[lb]{\smash{{\SetFigFont{9}{10.8}{\rmdefault}{\mddefault}{\updefault}$z_{q}^{-1} (D_{q}^{*}) = D (0, \sqrt{2 b} )^{*}$}}}}
\put(6526,-1936){\makebox(0,0)[lb]{\smash{{\SetFigFont{9}{10.8}{\rmdefault}{\mddefault}{\updefault}$z_{q}^{-1} ({\mathcal F}_{+})$}}}}
\put(10501,-2761){\makebox(0,0)[lb]{\smash{{\SetFigFont{9}{10.8}{\rmdefault}{\mddefault}{\updefault}$z_{q}^{-1} (\R \setminus 2 b \N)$}}}}
\put(1876,-1111){\makebox(0,0)[lb]{\smash{{\SetFigFont{9}{10.8}{\rmdefault}{\mddefault}{\updefault}$\C$}}}}
\put(10726,-4486){\makebox(0,0)[lb]{\smash{{\SetFigFont{7}{8.4}{\rmdefault}{\mddefault}{\updefault}$\hbox{Re} \, k$}}}}
\put(4651,164){\makebox(0,0)[lb]{\smash{{\SetFigFont{7}{8.4}{\rmdefault}{\mddefault}{\updefault}$\hbox{Im}\, k$}}}} 
\end{picture}%
\caption{{\bf Pre-images under $z_{q}$.}}
\end{center}
\end{figure}

\begin{prop}\label{lemA.1}
For each $N > 0$ the operator $(H_0-z)^{-1} : e^{-N \langle x_{3}
\rangle} L^2 (\R^{3}_{\bf x} ) \rightarrow e^{N \langle x_{3}
\rangle} L^2 (\R^{3}_{\bf x})$ has a holomorphic extension from
the open upper half plane to  ${\mathcal M}_{N}$. We denote
its holomorphic extension by $R_{0}(z)$.

Moreover, for $N>0$ and $\vp(\xp)=  \langle \xp \rangle^{- \alpha }$, $\alpha >1$, the holomorphic extension of
\begin{equation*}
T_{\vp}: z \mapsto \vp (\xp) e^{-N \langle x_{3} \rangle} (H_0-z)^{-1} e^{-N \langle x_{3} \rangle} ,
\end{equation*}
is holomorphic on ${\mathcal M}_{N}$ with values in the
Hilbert-Schmidt class $S_2$  on $L^{2} (\R^{3}_{\bf x})$.
\end{prop}

%
%
%

\begin{proof}
Since the kernel of $t_{N} (k^2)$ is given by
\begin{equation*}
e^{-N \langle x_{3} \rangle} \frac{i e^{i k |x_{3} -x_{3}' |}}{2 k } e^{-N \langle x_{3} ' \rangle} ,
\end{equation*}
the operator-valued function $k \mapsto t_{N} (k^2)$ has a holomorphic
extension from $\C_{+}$ to $\{ k \in \C^{*}; \ \Im k >-N \}$ in the
Hilbert-Schmidt class $S_2$ and in the trace class $S_{1}$ (see for instance
\cite{Fr97}). For $\Im k > 0$, we have the trace-class
estimate
\begin{align}
\| t_{N} (k^2) \|_1  &= \| e^{-N \langle x_{3} \rangle} (D_{x_3}
-k)^{-1}\|_2 \| (D_{x_3} + k)^{-1} e^{-N \langle x_{3} \rangle} \|_2 \leq
\frac{1}{2\pi N} \int_{\R} \frac{d \eta}{\eta^2 + |\Im k|^2}  \nonumber  \\
&= {\mathcal O}\left( |\Im k|^{-1}\right),
\label{eq2.6b}
\end{align}
and when moreover $\Re k^2<0$, we have the Hilbert-Schmidt estimate
\begin{align}
\| t_{N} (k^2) \|_2^2  = & \tr \Big( e^{-N \langle x_{3} \rangle} (D_{x_3}^2 -k^2)^{-1} (D_{x_3}^2 - \overline{k}^2)^{-1} e^{-N \langle x_{3} \rangle} \Big)  \nonumber  \\
\leq & \frac{1}{2\pi N} \int_{\R} \frac{d \eta}{(\eta^2 + |\Re k^2|)^2} = {\mathcal O}\left(|\Re k^2|^{-\frac32}\right), \label{eq2.6a}
\end{align}
where $\|\cdot\|_j$ stands for the norm in $S_j$, $j=1,2$.

By the definition of  ${\mathcal M}_{N}$, it follows that for any $q \in \N$,
the operator-valued  function $z \mapsto e^{-N \langle x_{3} \rangle}
(D_{x_3}^2 -z+2bq)^{-1} e^{-N \langle x_{3} \rangle} \in S_1$ can be
holomorphically extended from ${\mathcal F}_+$ to ${\mathcal M}_{N}$. We denote
its holomorphic extension to ${\mathcal M}_{N}$ also by $t_{N} (z-2b q)$.  Since $\{p_q, \: q \in \N\}$ is a family of orthogonal projectors,  we deduce  the holomorphic extension of (\ref{eq2.1}).

 Now, let us prove the existence of a holomorphic extension of $T_{\vp}$ in the
Hilbert-Schmidt class.
Let $z_0 \in {\mathcal M}_{N}$ be fixed, and $\Omega_0$ be a relatively compact
neighborhood of $z_0$. Since any path on ${\mathcal M}_{N}$ can enclose only a
finite number of Landau levels, there exists $q_0$  sufficiently large
(depending of $\Omega_0$) such that $q \geq q_0$ implies $t_{N} (z-2bq) =
e^{-N \langle x_{3} \rangle} (D_{x_3}^2 -z+2bq)^{-1} e^{-N \langle x_{3} \rangle}$.
Then for $q \geq q_0$ we have $\| t_{N} (z-2bq) \|  = {\mathcal O} \left( \langle q
\rangle^{-1} \right)$, and
furthermore, it follows from (\ref{eq2.6a}) that the identity
\begin{equation}\label{eq2.6}
\| t_{N} (z-2bq) \|_2 = \| e^{-N \langle x_{3} \rangle} (D_{x_3}^2
-z+2b q)^{-1} e^{-N \langle x_{3} \rangle} \|_2  = {\mathcal
  O} \Big( \langle q \rangle^{-\frac{3}{4}} \Big),
\end{equation}
holds for any $q \geq q_0$, uniformly with respect to $z \in \Omega_0$.

Next, we have
\begin{equation}  \label{eq2.7}
T_{\vp} (z) = \sum_{q =0 }^{q_0} \vp p_q \otimes t_{N} (z-2bq)+ \sum_{q > q_0} \vp p_q \otimes t_{N} (z-2bq),
\end{equation}
where $q_0$ is chosen as above. It is well known (see \cite{La}) that the orthogonal projection
$p_q$ admits an explicit integral kernel
\begin{equation} \label{ggg1}
{\mathcal P}_{q,b}(\xp,\xp^\prime)=\frac{b}{2\pi} L_q \left( \frac{b |\xp - \xp^\prime|^2}{2}\right)
\exp \Big( -\frac{b}{4} \big( |\xp - \xp^\prime|^2 + 2i(x_1x_2^\prime-x_1^\prime x_2) \big) \Big),
\end{equation}
where $L_q(t):=\frac{1}{q !}e^t \frac{d^q(t^qe^{-t})}{dt^q}$ are the Laguerre
polynomials. Note that ${\mathcal P}_{q,b}$ is constant onto the diagonal, i.e.
\begin{equation*}
{\mathcal P}_{q,b}(\xp,\xp)=\frac{b}{2\pi}, \quad \xp \in\R^2.
\end{equation*}
Further, if $U \in L^r(\R^2)$, $r \geq 1$, then $p_q U p_q$ is in the Schatten-von Neumann class $S_r$ (see Lemma 5.1 of \cite{Ra90}). In particular, $p_q \vp^2 p_q \in S_1$, and hence $\vp p_q \vp \in S_1$, and $\vp
p_q \in S_2$ with
\begin{equation*}
\| \vp p_q\|_2^2= \tr (\vp  p_q \vp) = \frac{b}{2 \pi} \int_{\R^2} \vp(\xp)^2 d\xp,
\end{equation*}
which is uniform with respect to $q$. Combining this with
(\ref{eq2.6}) and $p_{q} p_{k} = \delta_{q,k} p_{q}$, we deduce
that the infinite sum in \eqref{eq2.7} is convergent in $S_2$, and
hence $z \mapsto T_{\vp}(z) \in S_2$ has  a holomorphic extension
to ${\mathcal M}_{N}$. This  concludes the proof of Proposition
\ref{lemA.1}.
\end{proof}


For further references we formulate the following lemma which complements Proposition \ref{lemA.1}.

\begin{lem}  \label{lemdT}
 For $V$ satisfying \eqref{supexp} with $\mper >2$, the operator
\begin{equation*}
{\mathcal F}_+ \ni z \mapsto {\mathcal T}_V(z) :=  J |V|^{\frac12}(H_0-z)^{-1} |V|^{\frac12} \in S_2,
\end{equation*}
with $J: =  {\rm sign}\, V$ defined so that $J^2 = 1$,  has
 an analytic extension from ${\mathbb C}_+$ to ${\mathcal M}$,
denoted again by ${\mathcal T}_V(z)$. Moreover
the operator $d{\mathcal T}_V(z)/dz \in S_1$ is analytic on ${\mathcal M}$.
\end{lem}

\begin{proof}
The existence of the holomorphic extension in $S_2$ is a direct
consequence of Proposition \ref{lemA.1}, because for any $N>0$, we
have $|V|^{\frac12} = {\mathcal V} (\vp \otimes e^{-N \langle
x_{3} \rangle})$, with  a bounded operator ${\mathcal V}$.  In
order to prove $d{\mathcal T}_V(z)/dz \in S_1$, it suffices to
check that the series of general term $(\vp p_q  \vp) \otimes
dt_{N} (z-2bq)/dz$ converge in the trace class. Arguing as in the
proof of Proposition \ref{lemA.1}, we find that this convergence
follows from (\ref{eq2.6}), and the fact that for $z \in
\C\setminus [0, + \infty[$, we have $d t_{N} (z)/dz = e^{-N
\langle x_{3} \rangle} (D^2_{x_3}-z)^{-2} e^{-N \langle x_{3}
\rangle}$.
\end{proof}

\begin{rem}
The assumption $\mper>2$ could be weakened to
$\mper>1$ in the first part of the proof of the Lemma while it is
necessary for the second part.
\end{rem}

 Suppose \eqref{supexp} with $\mper > 0$. Using
\begin{equation*}
(H-z)^{-1} \big( 1 + V (H_{0} -z)^{-1} \big) = (H_0-z)^{-1},
\end{equation*}
we get
\begin{equation}
e^{-N \langle x_{3} \rangle} (H-z)^{-1} e^{-N \langle x_{3} \rangle} = e^{-N \langle x_{3} \rangle} (H_0-z)^{-1} e^{-N \langle x_{3} \rangle} \big( 1 + e^{N \langle x_{3} \rangle} V (H_{0} -z)^{-1} e^{-N \langle x_{3} \rangle} \big)^{-1} ,
\end{equation}
for $z \in {\mathcal F}_{+}$, $\im z \gg 1$. From Proposition \ref{lemA.1}, combined with \eqref{supexp}, the operator $V e^{N \langle x_{3} \rangle} (H_{0} -z)^{-1} e^{-N \langle x_{3} \rangle}$ is compact on $L^{2} (\R^{3}_ {\bf x})$. Then, the analytic Fredholm theorem proves the meromorphic extension of $\big( 1 + e^{N \langle x_{3} \rangle} V (H_{0} -z)^{-1} e^{-N \langle x_{3} \rangle} \big)^{-1}$ from ${\mathcal F}_{+}$ to ${\mathcal M}_{N}$. This now allows us to define the resonances of $H$.

\begin{prop}  \label{propA.1}
Suppose $V$ satisfies \eqref{supexp} with $\mper > 0$. Then  the
operator-valued function
\begin{equation*}
(H-z)^{-1} : e^{-N \langle x_{3} \rangle} L^2 (\R_{\bf x}^{3}) \rightarrow  e^{N \langle x_{3} \rangle} L^2 (\R_{\bf x}^{3}),
\end{equation*}
has a meromorphic extension from the open upper half plane to ${\mathcal
  M}_{N}$. Moreover, the poles and the  range of the residues of this extension do not depend on $N$.
\end{prop}

\begin{defi}
 We define the resonances of $H$ as the poles of the meromorphic extension of the resolvent $(H-z)^{-1}$, denoted by $R (z)$. The multiplicity of a resonance $z_0$ is defined by
\begin{equation}\label{eq2.10}
\mult (z_0): = \rank \frac{1}{2i\pi} \int_{\gamma} R (z) \, d z ,
\end{equation}
where $\gamma$ is a small positively oriented circle centered at
$z_0$.
\end{defi}

In the  sequel
we  will use also  the regularized
determinant ${\det}_2(I+A)$  defined for a Hilbert-Schmidt operator $A$ by
\begin{equation}\label{eq2.15}
{\det}_2(I+A):= \det \big( (I+A)e^{-A} \big).
\end{equation}

\begin{prop}\label{propA.2}
 Suppose $V$ satisfies \eqref{supexp} with $\mper >0$. The following assertions are equivalent:

{\rm (i)} $z_0 \in {\mathcal M}$ is a resonance of $H$,

{\rm (ii)} $z_0$ is a pole of $\vert V \vert^{\frac{1}{2}} R (z) \vert V \vert^{\frac{1}{2}}$,

{\rm (iii)} $-1$ is an eigenvalue of ${\mathcal T}_V(z_{0}) = J \vert V \vert^{\frac{1}{2}} R_{0} (z_{0}) \vert V \vert^{\frac{1}{2}}$.

\noindent Moreover, the rank of the residue of $\vert V
\vert^{\frac{1}{2}} R (z) \vert V \vert^{\frac{1}{2}}$ at
$z_{0}$ is equal to the multiplicity of the resonance of $H$.

Assume now that $V$ satisfies \eqref{supexp} with $\mper > 2$. Then
\begin{equation}  \label{rrp}
{\det}_2 \big( (H-z)(H_0-z)^{-1} \big) = {\det}_2 \big( I + {\mathcal T}_V(z) \big),
\end{equation}
has a analytic continuation from ${\mathcal F}_{+}$ to ${\mathcal M}$. Its zeros are the resonances of $H$, and if $z_0$ is a resonance, there exists a holomorphic function $f (z)$, for $z$ close to $z_0$, such that $f(z_0) \neq 0$ and
\begin{equation} \label{krein0}
\hbox{det}_2 \big( I + {\mathcal T}_V(z) \big) = (z-z_0)^{l(z_0)}f(z),
\end{equation}
with  $0 < l(z_0) = \mult (z_0)$ where $\mult (z_0)$ is the multiplicity of the resonance defined by \eqref{eq2.10}.
\end{prop}

\begin{rems}
(i) The main part of the proof of Proposition
\ref{propA.2} follows the arguments of \cite{SjBook}. To our best
knowledge the novelty is the proof of \eqref{rfv} (i.e. the
equality between the rank of the residue of $\vert V
\vert^{\frac{1}{2}} R (z) \vert V \vert^{\frac{1}{2}}$ at
$z_{0}$ and the multiplicity $\mult (z_0)$), and the equality  $l(z_0)
=\mult (z_0)$.

(ii) If $H$ and $H_0$ are two self-adjoint operators such that
$H-H_0 \in S_1$, the perturbation determinant ${\rm det}\left((H -
\lambda)(H_0-\lambda)^{-1}\right)$, ${\rm
  Im}\,\lambda > 0$, was introduced
by M. G. Krein in \cite{kr} (see also \cite[Section IV.3]{gk}). In the case
$H-H_0 \in S_r$ with $r>1$ the generalized perturbation determinant
 ${\rm det}_r\left((H - \lambda)(H_0-\lambda)^{-1}\right)$ was introduced in \cite{kop}.
In the last work, relatively Hilbert-Schmidt perturbations and the
corresponding generalized perturbation determinants ${\rm
  det}_2\left((H - \lambda)(H_0-\lambda)^{-1}\right)$ were considered as well; these determinants
are exactly of the same type as the one appearing in
(\ref{eq2.15}).

(iii) For  potentials $V$ compactly supported with respect to
$x_3$, the resonances can be defined as the poles of the
meromorphic extension of the resolvent:
\begin{equation}
(H-z)^{-1} :  L^2_{\rm comp} \big( \R_{x_3},L^2(\R^2_{\xp}) \big)
\rightarrow  L^2_{\rm loc} \big( \R_{x_3},L^2(\R^2_{\xp}) \big)
\end{equation}
from the open upper half plane to ${\mathcal M}$ (see \cite{SjBook} for the
Schr\"{o}dinger operator  $- \Delta + V$).
\end{rems}

\begin{proof}
Clearly, if $z_0$ is a pole of $\vert V \vert^{\frac{1}{2}}
R (z) \vert V \vert^{\frac{1}{2}}$ then it is a pole of
$R (z)$ and conversely according to the resolvent equation
\begin{align}
R (z) =& R_{0} (z) - R (z) V R_{0} (z) \nonumber  \\
=& R_{0} (z) - R_{0} (z) V R_{0} (z)  \nonumber  \\
&+ R_{0} (z) \vert V \vert^{\frac{1}{2}} J \vert V \vert^{\frac{1}{2}} R (z) \vert V \vert^{\frac{1}{2}} J \vert V \vert^{\frac{1}{2}} R_{0} (z) ,
\end{align}
if $z_0$ is a pole of $R (z)$ then it is a pole of $\vert V \vert^{\frac{1}{2}} R (z) \vert V \vert^{\frac{1}{2}}$, and (i) is equivalent to (ii).

From the resolvent equation we get
\begin{equation}  \label{pop}
\big( I + J \vert V \vert^{\frac{1}{2}} R_{0} (z) \vert V \vert^{\frac{1}{2}} \big) \big(I - J \vert V \vert^{\frac{1}{2}} R (z) \vert V \vert^{\frac{1}{2}} \big)=I.
\end{equation}
From Proposition \ref{propA.1} and Lemma \ref{lemdT}, we deduce that $z_0$ is a resonance if and only if $-1$ is an eigenvalue of ${\mathcal T}_V (z_0)$, and (ii) is equivalent to (iii).

Now we check the preservation of the multiplicity. Let $z_{0} \in {\mathcal M}$ be a resonance and $N>0$ be large enough to have $z_{0} \in {\mathcal M}_{N}$. For $z$ close to $z_0$, the resolvents, as operators from $L^2_{-N}:=e^{-N \langle x_{3} \rangle} L^2 (\R_{\bf x}^{3})$ to  $L^2_{N}$ , can be written
\begin{align}
R_{0}(z) &= \sum_{j \geq 0}M_j(z-z_0)^j  \\
R (z) &= (z-z_0)^{-L} A_{-L} + \cdots + (z-z_0)^{-1} A_{-1}+ {\hbox{Hol}}(z),\label{laurent}
\end{align}
were the last term is holomorphic in a neighborhood of $z=z_0$.
Classically, for $\gamma$ a small positively oriented circle
centered at $z_0$, we have
\begin{equation} \label{4}
A_{-j}= \frac{1}{2i\pi}\int_{\gamma} (z-z_0)^{j-1} R (z) \, d z,
\quad j \geq 1,
\end{equation}
$\mult (z_0)$ being the rank of $A_{-1}$, and
\begin{equation}  \label{eq2.25}
A_{-j} ( H - z_{0}) = ( H - z_{0}) A_{-j}=A_{-j-1}.
\end{equation}
%
%

Our next goal is to check the identities
\begin{equation} \label{gg1}
\rank (A_{-1})=\rank (V A_{-1} ),
\end{equation}
\begin{equation} \label{gg2}
\rank (A_{-1}^{*})=\rank (V A_{-1}^{*} ).
\end{equation}
Let us prove \eqref{gg1}. If this identity is false, there exists
a function $f$ such that $f=A_{-1} g$ and $V f=0$.

Since $f$ belongs to the range of $A_{-1}$, the distribution
 $H^{m} f$ is in $L_{\loc}^{2}(\R^{3})$, for any   $m \in \N$. In
particular, $H f = H_{0} f \in H^{2}_{\loc} (\R^{3})$, and hence
$f \in H^{4}_{\loc} (\R^{3}) \subset C^{2} (\R^{3})$.

Further, $V  f =0$ easily implies $V H f =0$. By recurrence, we obtain
\begin{equation}  \label{mlm}
V H^{n} f =0,
\end{equation}
for any $n \in \N$. Plugging \eqref{laurent} into the r.h.s. of the resolvent equation
\begin{equation*}
R (z) = R_{0} (z) - R_{0} (z) V R (z),
\end{equation*}
and integrating with respect to $z \in \gamma$, we find that \eqref{eq2.25} entails
\begin{equation} \label{ccc}
f=A_{-1} g= - \sum_{j=0}^{L-1} M_j V A_{-j-1} g= - \sum_{j=0}^{L-1} M_j V (H-z_0)^j f.
\end{equation}
Using \eqref{mlm} and \eqref{ccc}, we get $f=0$ which immediately yields
\eqref{gg1}. Identity \eqref{gg2} can be proved exactly in the same way.

Applying \eqref{gg1} -- \eqref{gg2}, we get $\rank (A_{-1})=\rank ( A_{-1} V )$,
next ${\rank}\,(A_{-1})= {\rank}\,( A_{-1} V) = {\rank}\, (A_{-1} \vert
V \vert^{\frac{1}{2}})$, and, moreover, find that $\vert V \vert^{\frac{1}{2}}$ is
injective on the  range of $ (A_{-1})$. Thus we obtain
\begin{equation}  \label{rfv}
\rank (A_{-1}) = \rank \big( \vert V \vert^{\frac{1}{2}} A_{-1} \vert V \vert^{\frac{1}{2}} \big),
\end{equation}
which implies that the multiplicities agree.

We now prove the second part of the proposition. Let us recall
that, if $A$ is a bounded operator and if $B$ is a trace class
operator on some separable Hilbert space, we have $\det (I +AB) =
\det (I+BA)$. Moreover, for $A$ bounded and $B$ Hilbert-Schmidt,
we have
\begin{equation} \label{vvb}
{\det}_{2} (I+AB) = {\det}_{2} (I+BA).
\end{equation}
Writing, for  $z \in {\mathcal F}_{+}$,
\begin{equation*}
(H-z)(H_0-z)^{-1} = I + V (H_0-z)^{-1},
\end{equation*}
where $J \vert V \vert^{\frac{1}{2}} (H_0-z)^{-1}$ is holomorphic
on ${\mathcal F}_{+}$, with value in the Hilbert-Schmidt class, and
 using \eqref{vvb}, we get
\begin{equation*}
{\det}_2 \big( (H-z)(H_0-z)^{-1} \big) = {\det}_2 \big( I + \vert V \vert^{\frac{1}{2}} J \vert V \vert^{\frac{1}{2}} (H_0-z)^{-1} \big) = {\det}_2 \big( I + {\mathcal T}_V(z) \big) .
\end{equation*}
From Lemma \ref{lemdT}, this determinant has an analytic extension from ${\mathcal F}_{+}$ to ${\mathcal M}$ and vanishes if and only if $z_{0}$ is a resonance of $H$.
Then, there exists a holomorphic function $f (z)$, for $z$ close to $z_0$, such that $f(z_0) \neq 0$ and
$$\hbox{det}_2 \big( I + {\mathcal T}_V(z) \big) = (z-z_0)^{l(z_0)}f(z_0).$$
In order to prove that  $l(z_0) = \mult (z_0)$ where $\mult (z_0)$ is the
multiplicity of the resonance defined by \eqref{eq2.10}, we need
the following

\begin{lem}\label{lempia=a}
The operator
$$\Pi_{-1} =- A_{-1} \sum_{L-1 \geq j,k\geq 0} (H-z_0)^j V M_{j+k+1}V(H-z_0)^k,$$
is well defined in ${\mathcal L}({{\mathcal H}}^{L-1}_N)$
for any $N$ where ${{\mathcal H}}^L_N$ is the Hilbert space
$${{\mathcal H}}^L_N:= \{u \in L^2_N(\R^3)=L^2(\R^3, e^{-N|x_3|}dx) \; {\hbox{ such that }} H^ku\in  L^2_N(\R^3), \; \forall k \leq L\},$$
equiped with the norm $\sum_{0\leq k \leq L}\| H^k u \|_{L^2_N}$.

For any $k \in \N$, in ${\mathcal L}(L^2_{-N},{{\mathcal H}}^k_N)$ we have
\begin{equation}\label{eqpiaa}
\Pi_{-1} A_{-1} = A_{-1}.
\end{equation}
Here ${\mathcal L} (A,B)$ (resp. ${\mathcal L} (A)$) denotes the space
of linear bounded operator from $A$ to $B$ (resp. $A$).
\end{lem}

\begin{proof}[Proof of Lemma \ref{lempia=a}]

We recall that from \eqref{ccc},
\begin{equation}\label{eqlm0}
A_{-1} = - \sum_{j\geq 0} M_j V A_{-j-1},
\end{equation}
with  the convention that $A_{-j}=0$ for $j>L$. On the other hand,
the resolvent equation
$$R_{0} (z) = R (z) \big( I + V R_{0} (z) \big),$$
yields
\begin{equation}\label{eqtmj}
M_j = \sum_{k\leq j} A_k {\widetilde M}_{j-k}
\end{equation}
with ${\widetilde M}_0=I + V M_0$,  ${\widetilde M}_{j}=V M_j$  for $j\geq
1$, and the equation $R_{0} (z) = \big( I +
R_{0} (z) V \big) R (z)$ implies for any  $k \geq 0$:
\begin{equation}\label{eqttmj}
0 = \sum_{ j\geq k} {\widetilde{\widetilde{M}}}_{j-k} A_{-j-1}
\end{equation}
with  ${\widetilde{\widetilde{M}}}_{0}=I+M_0V$, $
{\widetilde{\widetilde{M}}}_{j}=M_jV$, $ {\widetilde{\widetilde{M}}}_{-j}=0$  for
$j\geq 1$. In the above equality we use again the convention that
$A_{-j}=0$ for $j>L$. By inserting (\ref{eqtmj}) into
(\ref{eqlm0}), we deduce
$$A_{-1} = - \sum_{j\geq 0}\sum_{k\leq j} A_k {\widetilde M}_{j-k}  V A_{-j-1}.$$
Since ${\widetilde M}_{j}V=V{\widetilde{\widetilde{M}}}_{j}$ and
$A_{-1}(H-z_0)^j = A_{-1-j}$, relation \eqref{eqttmj} implies
$$A_{-1} = - \sum_{j\geq 0}\sum_{k<0} A_k V M_{j-k}  V A_{-j-1}= \Pi_{-1}A_{-1}.$$
This concludes the proof of Lemma \ref{lempia=a}.
\end{proof}

Let us now complete the proof of Proposition \ref{propA.2}. It
follows from Lemma \ref{lempia=a} that $\rank \Pi_{-1} =\rank
A_{-1}$ and (\ref{eqpiaa}) implies that
 $\Pi_{-1}\Pi_{-1}=\Pi_{-1}$. Consequently, we have
\begin{equation}\label{mm0}
\mult (z_0)= \rank A_{-1}=\tr \; \Pi_{-1}.
\end{equation}
On the other hand, by the definition of $l(z_0)$, we have
\begin{equation}\label{lz0}
l(z_0)= \frac{1}{2i\pi} \int_{\gamma} \partial_z \ln {\det}_2\Big(1 + {\mathcal T}_V (z)\Big) dz.
\end{equation}
Further, we have
$$ \partial_z \ln {\det}\Big(1 + T(z)\Big)= \tr \Big( (1+T(z))^{-1}\partial_z T(z) \Big), \quad z \in \Omega,$$
for any operator-valued holomorphic function $\Omega \ni z \mapsto
T(z) \in S_1$. Therefore,
$$ \partial_z \ln {\det}_2\Big(1 + {\mathcal T}_V (z)\Big)=\tr \Big( (1+{\mathcal T}_V(z))^{-1}\partial_z {\mathcal T}_V(z) \Big)- \tr  \Big( \partial_z {\mathcal T}_V(z) \Big).$$
According to Lemma \ref{lemdT}, $\partial_z {\mathcal T}_V(z)$ is holomorphic in the trace class, then its integral on $\gamma$ vanishes and (\ref{pop}) yields:
$$l(z_0)= - \frac{1}{2i\pi} \int_{\gamma} \tr\Big( J |V|^{\frac12} R (z) V\partial_z R_0(z) |V|^{\frac12}\Big) dz.$$
By definition of $A_{-k}$ and $M_j$,
we obtain:
$$l(z_0)= - \tr \Big(  \sum_{ L \geq k\geq 1} J |V|^{\frac12} A_{-k} k V  M_k|V|^{\frac12}\Big),$$
where the trace is on ${\mathcal L}(L^2)$. Thanks to
(\ref{eq2.25}), we have $A_{-k} = A_{-1} (H-z_0)^{k-1}$ in
${\mathcal L}(L^2_{-N}, L^2_N)$ and using the cyclicity of the
trace,
\begin{equation}\label{eqlz0}
l(z_0)= - \tr \Big(  \sum_{ L \geq k\geq 1}  A_{-1} (H-z_0)^{k-1}k V  M_k V\Big).
\end{equation}
Here the trace is in ${\mathcal L}( L^2_N)$, but since the range
of $A_{-1}$ is in any ${{\mathcal H}}^j_N$, $j \geq 1$, the last trace is
also in any ${{\mathcal H}}^j_N$. At last, combining the cyclicity of the
trace, with (\ref{eq2.25}), (\ref{mm0}) and (\ref{eqlz0})  we
deduce
\begin{equation*}
\mult (z_0) = \tr ( \Pi_{-1})= \tr \Big(- \sum_{0 \leq j,k \leq L-1} A_{-1} (H-z_0)^{k+j}V M_{k+j+1}V \Big)=l(z_0).
\end{equation*}
\end{proof}

\section{Resonances near  the Landau levels}

In this section we assume that $V$ satisfies \eqref{supexp} with
$\mper >2$, and study  the resonances localized in $D_{q}^{*}$, the
neighborhood of the Landau level $2b q$ introduced in Definition
\ref{a36}. Recall that $D^{*}_{q}$ can be parametrized by $z_{q}
(k)$ defined in \eqref{a35}.

According to the previous section, these resonances can be identified with
the points $z$ where the determinant $\det_2(I+ T_V(z))$ vanishes. Note that ${\mathcal T}_V(z)$ is the holomorphic extension of
\begin{equation}\label{eq2.10b}
J |V|^{\frac12}(H_0-z)^{-1} |V|^{\frac12}= \sum_{j \in \N}  J
|V|^{\frac12}  P_j (H_0-z)^{-1} |V|^{\frac12}, \quad z \in
{\mathcal F}_+,
\end{equation}
where $P_j= p_j \otimes I_{x_3}$, $j \in {\mathbb N}$.


In order to study  the resonances near a Landau level $2bq$ we split
${\mathcal T}_V(z)$ into two parts:
\begin{equation*}
{\mathcal T}_V(z) =  J |V|^{\frac12}  P_q R_{0} (z) |V|^{\frac12} +  \sum_{j \neq q}
J |V|^{\frac12}  P_j R_{0} (z) |V|^{\frac12}.
\end{equation*}
By Proposition \ref{lemA.1}, the second term  in the r.h.s. is
holomorphic in a neighborhood of $2bq$ with values in $S_2$. Let
us consider the first term for $z=z_{q} (k)$.  The series
expansion with respect to $k$ of the kernel of the operator
$t_{N}$ (see \eqref{mml} and \eqref{tchi}) allows  us to write
$t_{N}$  as the sum
\begin{equation}\label{eq3.10}
t_{N}(k^2)=\frac{1}{k} t_1 + r_1(k),
\end{equation}
where $t_1 : L^2(\R) \to L^2(\R)$ is the rank-one operator defined by
\begin{equation}\label{deft1}
t_1 u : = \frac{i}{2}
\big\langle u , e^{-N \langle . \rangle} \big\rangle e^{-N \langle x_{3} \rangle},
\end{equation}
and $r_1(k)$ is the Hilbert-Schmidt operator with  integral kernel
\begin{equation}\label{defr1}
{\mathcal R}_1(x_3,x_3')= e^{-N \langle x_{3} \rangle} \;i \frac{ e^{i k |x_3-x_3'|}-1}{2k}
e^{-N \langle x_{3}' \rangle} .
\end{equation}
 Clearly, the operator-valued function ${\mathbb C} \ni k \mapsto r_1(k)
\in S_2$ is analytic. Putting together the above considerations, we obtain the
following

\begin{prop}\label{propB.2}
Suppose  $V$  satisfies \eqref{supexp} with $\mper >2$. For $k \in \C^*$, $|k|< \sqrt{2b}$, we have:
\begin{equation}\label{eq3.5}
{\mathcal T}_V( z_{q} (k)) = \frac{iJ}{k} B_q + A(k),
\end{equation}
where $J = {\rm sign}\,V$, $B_q$ is the positive self-adjoint
operator
\begin{equation}
B_{q} = \frac{1}{2} \vert V \vert^{1/2} P_{q} \vert V \vert^{1/2},
\end{equation}
and $A(k) \in S_2$ is the holomorphic operator defined on
$\{k \in \C, \; |k|< \sqrt{2b}\}$ by 
\begin{equation}\label{defAk}
A(k) = JA_q(k) + J \sum_{j \neq q}   |V|^{\frac12}  P_j R_0(z_q(k))|V|^{\frac12},
\end{equation}
where $A_q(k)$ is the operator with  integral kernel
\begin{equation*}
{\mathcal K}_{A_q}(\xp,x_3;\xp',x_3')= |V(\xp,x_3)|^{\frac12}
{\mathcal P}_{q,b} (\xp,\xp') \frac{1- e^{i k |x_3-x_3'|}}{2ik}
|V(\xp',x_3')|^{\frac12}.
\end{equation*}
Here, ${\mathcal P}_{q,b}$ is the integral kernel of the orthogonal
projection $p_q$ written in \eqref{ggg1}.
\end{prop}

 Since there exists an operator $C: L^2(\R^3) \rightarrow
L^2(\R^2)$ such that $B_q=C^*C$ and $CC^*= \frac12 p_q W p_q$ with
\begin{equation}\label{defW}
W(\xp)=\int_{\R} |V(\xp,x_3)| dx_3,
\end{equation}
(see  \cite{sob1} for $q=0$, and the proof of Proposition 5.3
of \cite{FR}  for any $q \in {\mathbb N}$), then for any
$s>0$ we have
\begin{equation}\label{nBnW}
n_+(s; B_q) = n_+(2s;p_q W p_q),
\end{equation}
where for a compact self-adjoint operator $A$,  we set $n_+(s;A) = \rank \,
 {\bf 1}_{(s,+\infty)}(A)$.

\begin{rem}
Using (\ref{pop}) and (\ref{eq3.5}), we can prove that each Landau
level is an essential singularity of the resolvent of $H$, but it
is not sufficient to deduce the existence of an infinite number of
resonances near  the Landau levels.
\end{rem}

  In the case where the decay of $U$ at infinity is
regular enough, the asymptotic distribution of the eigenvalues of
Toeplitz-type operators $p_q U p_q$ is well known. The following
three lemmas describe the eigenvalue asymptotics for $p_q U p_q$
in the case of power-like decay, exponential decay, or compact
support of $U$, respectively.

\begin{lem}\label{lemRa90} {\rm (Theorem 2.6 of \cite{Ra90})}
Let the function $U \in C^1(\R^2)$ satisfy the estimates
$$0 \leq U(\xp) \leq C_1 \langle \xp \rangle^{-\alpha}, \qquad |
\nabla U(\xp) | \leq C_1 \langle \xp \rangle^{-\alpha-1}, \quad
\xp \in \R^2,$$ for some $\alpha >0$ and $C_1>0$. Assume, moreover
that
$$U(\xp) = u_0(\xp/|\xp|) \; |\xp|^{-\alpha} (1 + o(1)), \quad
|\xp| \rightarrow \infty,$$ where $u_0$ is a continuous function
on $S^1$ which is non-negative and does not vanish identically.
 Then for each $q \in \N$ we have
$$n_+(s;p_q U p_q) = C_{\alpha} s^{-2/\alpha}  (1 + o(1)), \quad s \searrow 0,
$$
where
\begin{equation} \label{ggg12}
C_{\alpha} : = \frac{b}{4 \pi} \int_{S^1}
u_0(t)^{2/\alpha} dt.
\end{equation}
\end{lem}

\begin{lem}\label{lem1RW} {\rm (Theorem 2.1 of \cite{RW})}
Let $0\leq U \in L^\infty(\R^2)$. Assume that
$$\ln U(\xp) = -\mu |\xp|^{2\beta}(1 + o(1)), \quad |\xp|
\rightarrow \infty,$$ for some $\beta >0$, $\mu >0$. Then for each
$q \in \N$ we have
$$n_+(s;p_q U p_q) = \varphi_\beta(s) (1+o(1)), \quad s \searrow 0,$$
where
$$\varphi_\beta(s):= \left\{ \begin{array}{ccc}
                            \frac{b}{2 }\mu^{-\frac{1}{\beta}} |\ln
                            s|^{\frac{1}{\beta}} & {\rm if} &
                            0<\beta<1,\\
                            \frac{1}{\ln (1+ 2 \mu/b)} |\ln s|& {\rm if} &
                            \beta=1,\\
                            \frac{\beta}{\beta-1}(\ln|\ln s|)^{-1} |\ln
                            s| & {\rm if} &
                            \beta>1,\\
                            \end{array}
                    \right. \qquad 0< s <  e^{-1}.$$

\end{lem}

\begin{lem}\label{lem2RW} {\rm (Theorem 2.4 of \cite{RW})}
Let $0\leq U \in L^\infty(\R^2)$. Assume that the support of $U$
is compact and there exists a constant $C>0$ such that $U\geq C$
on an non-empty open subset of $\R^2$. Then for each $q \in \N$ we
have
$$n_+(s;p_q U p_q) = \varphi_\infty(s) (1+o(1)), \quad  s \searrow 0,$$
where
$$\varphi_\infty(s):= (\ln|\ln s|)^{-1} |\ln
                            s|, \quad 0< s <  e^{-1}.$$
\end{lem}

\begin{rem}
In the recent preprint \cite{FiPu} a sharper version of the result
of Lemma \ref{lem2RW} has been obtained, containing three
asymptotic terms  as $s \searrow 0$ of $n_+(s;p_q U p_q)$ provided
that the support of $U$ is compact, and $U$ satisfies some
additional technical assumptions. In particular, the asymptotic
expansion of $n_+(s;p_q U p_q)$ obtained in  \cite{FiPu} recovers
{\em the  logarithmic capacity} of the support of $U$.
\end{rem}

The above lemmas imply some useful properties of $B_q$
summarized in the following

\begin{cor}\label{lemBq} Let $s>0$. For $V$ satisfying (\ref{eq1.1})
  with $\mper >2$, the operator $B_q$ is a trace class operator  with  $n_+(s,B_q)= {\mathcal O}(s^{-2/\mper})$ for $s>0$ small enough.
 For $j\in \N^*: = \{1,2,\ldots\}$, the operator-valued
 functions
\begin{equation}\label{ggg7}
 \C\setminus (\mp i [0, + \infty[)\, \ni k \mapsto
{\mathcal B}(k) = {\mathcal B}_{q,j}^{\pm}(k) : = \frac{iB_q}{k} \Big(I\pm
    \frac{iB_q}{k}\Big)^{-j} \in S_1
\end{equation}
are holomorphic.  Their Hilbert-Schmidt  norms (p=2) and
trace-class norms (p=1) satisfy
    the estimates
\begin{equation}\label{eqesti}
\|{\mathcal B}(k)\|_p \leq  c(\theta)^j \; \sigma_p(|k|)^{\frac{1}{p}},
\end{equation}
where $\theta = {\rm Arg}\, k$, $ c(\theta) = (1 - (\sin
\theta)_-)^{-\frac12}$, $u_-:= \max\{-u, 0\}$ if $u \in \R$, and
\begin{equation}\label{defsigmap}
 \sigma_p(s):= \left\| \frac{B_q}{s} \Big(I +\frac{B_q^2}{s^2}\Big)^{-1/2}\right\|_p^p= {\mathcal O}(s^{-2/\mper}), \quad s>0 .
\end{equation}

Further,  for $s>0$, $p \geq 1$,  we have
\begin{equation}\label{tnequivsigma}
2^{-p/2} \,  \widetilde n_p(s) \,  \leq \,  \sigma_p(s) \,  \leq \, \widetilde n_p(s) + n_+(s, B_q),
\end{equation}
where
\begin{equation}\label{deftn}
 \widetilde n_p(s):= \left\| \frac{B_q}{s} {\bf 1}_{[0,s]}(B_q) \right\|_p^p,
 \quad s>0, \quad p \geq 1.
\end{equation}

 Moreover, for $W$ defined by (\ref{defW})  satisfying  the
assumptions of Lemma \ref{lemRa90} with $\alpha > 2$, the estimates
\begin{equation}\label{asymsigma}
\sigma_p(s) = C_{\alpha,p}s^{-\frac{2}{\alpha}}(1 + o(1)), \quad  \widetilde n_p(s)= \widetilde C_{\alpha,p}s^{-\frac{2}{\alpha}}(1 + o(1)) \quad  s \searrow
0,
\end{equation}
hold with some $C_{\alpha,p} > 0$,  $\widetilde C_{\alpha,p} > 0$,
$p=1,2$. Finally, if the assumptions of Lemma \ref{lem1RW} or of
Lemma \ref{lem2RW} hold for $W=U$, we have
\begin{equation}\label{asymsigma1}
\sigma_p(s) = \varphi_\beta(s)(1 + o(1)), \quad  \widetilde n_p(s) = o(\varphi_\beta(s)) \quad s \searrow 0,
\end{equation}
the functions $\varphi_\beta(s)$, $0< \beta \leq \infty$, being defined in Lemma \ref{lem1RW} or in  Lemma \ref{lem2RW}.

\end{cor}

\begin{proof}
 By (\ref{nBnW}), $B_q \in S_1$ if and only if $p_q W p_q
\in S_1$. For $V$ satisfying (\ref{eq1.1}), we have $0 \leq W(\xp)
\leq C \langle \xp \rangle^{-\mper}$.  Therefore,  $W \in
L^1(\R^2)$, and hence $p_q W p_q \in S_1$, and then $B_q \in S_1$.  According to  Lemma \ref{lemRa90}, $n_+(s, p_q C \langle \xp \rangle^{-\mper}p_q)$ behaves like $s^{-2/\mper}$ as $s \searrow 0$, then  $n_+(s,B_q)= {\mathcal O}(s^{-2/\mper})$ for $s>0$ small enough.
Taking also into account that $B_q \geq 0$, we conclude that the
operator-valued functions ${\mathcal B}$ defined in \eqref{ggg7} are
holomorphic.  Let us now estimate  their norms in $S_1$
and in $S_2$. For $k = |k| e^{i \theta}$, we have
$${\mathcal B}^*{\mathcal B}= \frac{B_q^2}{|k|^2} \Big(I+ \frac{B_q^2}{|k|^2}  \pm 2\sin \theta \frac{B_q}{|k|}\Big)^{-j}.$$
 Next,
\begin{equation} \label{ggg2}
\|{\mathcal B}\|_{S_p}^p=
- \int_0^\infty f_{j,p}(\frac{u}{|k|},\pm\theta) d n_+(u;B_q) =
- \int_0^\infty f_{j,p}(s,\pm\theta) d n_+(s;\frac{B_q}{|k|}),
\end{equation}
where $f_{j,p}(u,\theta):=u^p(1+u^2+ 2u\sin \theta )^{-jp/2}$.
Evidently, for $\theta \neq  - \pi/2$ and  $u\geq 0$ we have
\begin{equation} \label{ggg3}
f_{j,p}(u,\theta)\leq  c(\theta)^{jp} f_p(u)
\end{equation}
where $f_p(u):= u^p(1+ u^2)^{-p/2}$, $p=1, 2$. Finally,
\begin{equation}\label{eq3.14}
\sigma_p(s)=\tr \Big( \frac{B_q^p}{s^p} \Big(I +\frac{B_q^2}{s^2}\Big)^{-p/2}
\Big) = - \int_0^\infty f_p(u) d n_+(u;\frac{B_q}{s}).
\end{equation}
Now the combination of \eqref{ggg2}, \eqref{ggg3}, and
\eqref{eq3.14}, yields \eqref{eqesti}.
 We have also
\begin{equation}\label{eq3.14b}
\widetilde n_p (s) = - \int_0^1 u^p d n_+(u;\frac{B_q}{s}).\end{equation}
Then \eqref{tnequivsigma} is a consequence of the  elementary inequalities
$$2^{-p/2} u^p {\bf 1}_{[0,1]}(u) \, \leq f_p(u) \, \leq \,  u^p {\bf
  1}_{[0,1]}(u)  + {\bf 1}_{]1,+\infty[}(u).$$

  In order to prove \eqref{asymsigma} - \eqref{asymsigma1}, we first note that since $\lim_{u \downarrow 0}\, u\, n_+(u,B_q)=0$,  relations \eqref{eq3.14},  \eqref{eq3.14b} and (\ref{nBnW}) imply
\begin{equation}\label{ggg4}
\sigma_p(s)  =  \int_0^\infty f_p'(u) n_+(2su;p_qWp_q) du,
\end{equation}
\begin{equation}\label{ggg4b}
 \widetilde n_p (s) = \int_0^1 p u^{p-1} ( n_+(2su;p_qWp_q)- n_+(2s;p_qWp_q))du .
\end{equation}
Then for $W$  satisfying the
assumptions of Lemma \ref{lemRa90} or of Lemma \ref{lem1RW} or of
Lemma \ref{lem2RW} we deduce   the asymptotic properties claimed.
\end{proof}

\begin{prop}\label{propB.3}
Suppose that $V$ satisfies (\ref{supexp}) with $\mper >2$. For $0<s<|k|<s_0$ with $s_0$ sufficiently
small,  $z_{q} (k) \in D_q^*$ is a resonance of $H$ if and only if $k$ is a
zero of
\begin{equation}\label{defDk}
D(k,s)=\det \Big( I + K(k,s)  \Big),
\end{equation}
where $K(k,s)$ is a  finite-rank operator satisfying
$$
\rank K(k,s)= {\mathcal O}\Big(n_+(s;p_q W p_q)+1 \Big),
\quad \| K(k,s)
\|={\mathcal O}(s^{-1}),
$$
uniformly with respect to $s<|k|<s_0$.

 Moreover, for $\Im k^2 > \delta >0$, the operator $I + K(k,s)$ is invertible
with
$$\|(I + K(k,s))^{-1}\|= {\mathcal O}(\delta^{-1}),$$
uniformly with respect to $s<|k|<s_0$, $\Im k^2 > \delta$.
\end{prop}

\begin{proof}
 By Proposition \ref{propA.2} --
\ref{propB.2}, for $s<|k|\leq s_0 < \sqrt{2b}$,  $z_q(k)$ is a
resonance of $H$ if and only if $k$ is a zero of $\det_2
(I+\frac{iJ}{k} B_q + A(k))$.

Since $k \mapsto A(k)$ is holomorphic near $k=0$ with value in
$S_2$, for $s_0$ sufficiently small, there exist a
finite-rank operator $A_0$ independent of $k$ and $\widetilde{A}(k)$
 holomorphic  near $k=0$  in $S_2$ with $\| \widetilde{A}(k)\|\leq  \frac14$, $|k|\leq s_0$ such that
$$A(k)=A_0 + \widetilde{A}(k).$$

 Further, let us decompose the self-adjoint positive
operator $B_q$ into a trace-class
 operator whose norm is
bounded by $s/2$, and an operator of rank $n_+(s/2;B_q)$, namely
\begin{equation}\label{eq3.15}
B_q = B_q  {\bf{1}}_{[0, s/2]}(B_q ) + B_q {\bf{1}}_{]s/2, +
\infty[} (B_q) .
\end{equation}

Since $\|  \frac{iJ}{k} B_q  {\bf{1}}_{[0, s/2]}(B_q ) +
\widetilde{A}(k) \| \leq \frac34$, for $0<s<|k|<s_0$, we have 
$$\det \Big( (I+ \frac{iJ}{k} B_q  {\bf{1}}_{[0, s/2]}(B_q ) +
\widetilde{A}(k))  \; e^{-{\mathcal T}_V(z_q(k))} \Big) \neq 0.$$

It follows that for $0<s<|k|<s_0$, the zeros of  $ \det_2 \Big( I +
{\mathcal T}_V(z_q(k)) \Big)$  are the zeros of $D(k,s)$ defined by
(\ref{defDk}) with
$$ K(k,s)=\Big(\frac{iJ}{k} B_q {\bf{1}}_{]s/2, +
\infty[} (B_q) + A_0 \Big) \Big(I+\frac{iJ}{k} B_q  {\bf{1}}_{[0,
s/2]}(B_q ) + \widetilde{A}(k)\Big)^{-1}.$$ The rank of this operator
is bounded by ${\mathcal O}(n_+(s/2;B_q)+1) = {\mathcal
O}(n_+(s;p_qWp_q)+1)$ (see (\ref{nBnW})) and its norm is bounded
by ${\mathcal O}(|k|^{-1})$.

 At last, by the definition of ${\mathcal T}_V(z)$ and of
$K(k,s)$, we have 
$$I+K(k,s) = \Big( I +{\mathcal T}_V(z_q(k))\Big) \Big(I+\frac{iJ}{k} B_q  {\bf{1}}_{[0,
s/2]}(B_q ) + \widetilde{A}(k)\Big)^{-1},$$
provided that
$0<s<|k|<s_0$.
 By the
resolvent equation  (\ref{pop}), the operator  $ I+ {\mathcal T}_V(z)
$ is invertible for $\Im
z>\delta$, and
$$  \Big( I + {\mathcal T}_V(z) \Big)^{-1}=  I -J |V|^{\frac12}(H-z)^{-1}|V|^{\frac12} .$$

Then  $I+K(k,s)$ is
invertible for $\Im k^2> \delta$, $0<s<|k|<s_0$, and
$$\| ( I + K(k,s) )^{-1} \| = {\mathcal O}(1 + \| |V|^{\frac12}(H-z_q(k))^{-1}|V|^{\frac12}\|)=
{\mathcal O}(1 + |\Im k^2|^{-1})$$ which concludes the proof of
Proposition \ref{propB.3}.
\end{proof}

By the properties of $K(k,s)$ (see Proposition \ref{propB.3}) for
$0<s<|k|<s_0$, we have:

\begin{equation}\label{ODk}
D(k,s)=  \prod_{j=1}^{{\mathcal O}(n_+(s;p_qWp_q)+1)}\;
(1+\lambda_j(k,s)) = {\mathcal O}(1) \exp \Big( {\mathcal
O}(n_+(s;p_qWp_q)+1) | \ln s|\Big),
\end{equation}
uniformly with respect to $(k,s)$, where $\lambda_j(k,s)$ are the
eigenvalues of $K(k,s)$ which  satisfy $\lambda_j(k,s)=  {\mathcal O}(|s|^{-1})$.

Moreover, since
$$ D(k,s)^{-1} =\det \Big(( I +K)^{-1} \Big)= \det \Big( I -K ( I + K
)^{-1}\Big),$$ for $\Im k^2> \delta >0$, and for $0<s<|k|<s_0$, we
have
\begin{equation}\label{-ODk}
| D(k,s) | \geq C \exp \Big(  -C (n_+(s;p_qWp_q)+1) (| \ln \delta|+ | \ln s|)\Big),
\end{equation}
uniformly with respect to $(k,s)$.

The following lemma contains a version of the well known Jensen inequality
which is suitable for our purposes.

\begin{lem}\label{lemJensen}
Let $\Omega$ be a simply connected
sub-domain of $\C$ and
let $g$ be a holomorphic function in $\Omega$ with continuous extension to
$\overline{\Omega}$. Assume there exists $z_0 \in \Omega$ such
that $g(z_0) \neq 0$ and $g(z) \neq 0$ for $z \in \partial
\Omega$.  Let $z_1, z_2, ...,z_N \in \Omega$ be the zeros of $g$
repeated according to their multiplicity. For any domain $\Omega'
\subset \subset \Omega$, there exists $C'>0$ such that
$N(\Omega',g)$, the number of zeros $z_j$ of $g$ contained in
$\Omega'$, satisfies
$$N(\Omega',g) \leq C \Big( \int_{\partial \Omega} \ln |g(z)| d z - \ln |g(z_0)|\Big).$$
\end{lem}
\begin{proof}
 First, let us recall the classical Jensen inequality
$$N(B(0,\nu R),G) |\ln \nu | \leq \frac{1}{2 \pi} \int_0^{2
\pi} \ln |G(R e^{i \theta})| d \theta - \ln |G(0)|,$$ valid for any $0<
\nu < 1$, and
for a function $g=G$ satisfying the assumptions of
the lemma in $\Omega =B(0,R):=\{z \in \C; \; |z| <R \}$, and
$z_0=0$.

 Now, let $f: \; B(0,R) \rightarrow \Omega$ be a
bijective analytic function such that $f(0)=z_0$ and $f(\partial
B(0,R)) = \partial \Omega$. Then $G=gof$ satisfies the assumptions of
the lemma in $\Omega =B(0,R)$, with $z_0=0$, and we have the above
formula. Since $f$ is a bijection and $f(\partial B(0,R)) =
\partial \Omega$, for $\Omega' \subset \subset \Omega$ there
exists $0<\nu<1$ such that $\Omega' \subset f( B(0,\nu R))$, which implies
the claim of the lemma.
\end{proof}

Applying this lemma to the function  $g(k):= D(r\, k,r
)$,   on $\Omega:=\{ k \in \C; \; 1<|k|<2 , \; \frac{\pi}{3}< \hbox{Arg } k<2 \pi +\frac{\pi}{6}\}$
 with $\Im k_0^2 > \delta> 0 $,
we deduce from (\ref{ODk}), (\ref{-ODk}) the following upper bound
on the number of resonances near the Landau levels.

\begin{thm}\label{thmUB} {\bf Upper bound.}
Suppose that $V$ satisfies (\ref{supexp}) with $\mper >2$.   Then there exists
$r_0>0$, such that for any  $0<r<r_0$,
$$\# \{ z = z_{q} (k) \in \res \cap D_{q}^{*}; \ r<\vert  k \vert < 2 r \} =
{\mathcal O}(n_+(r, p_q W p_q) | \ln r|),$$ where $W$ is given by
(\ref{defW}), and $n_+(s; p_q W p_q)$ is the counting function
satisfying asymptotic relations depending on the decay of $W$,
described in Lemmas \ref{lemRa90}, \ref{lem1RW},  and
\ref{lem2RW}. In particular, under our assumptions we have always
$n_+(s; p_q W p_q)= {\mathcal O}(s^{-2/\mper})$, and for $V$
compactly supported, we have $n_+(s; p_q W p_q)= {\mathcal
O}((\ln|\ln s|)^{-1}| \,\ln s|)$.
\end{thm}

\begin{rem}
Instead of the three-dimensional case considered in the present
paper, it is possible to consider a general $n$-dimensional
Schr\"odinger operator with non vanishing constant magnetic field
$B$ which can be regarded as a real antisymmetric matrix acting in
$\R^n$. Set $2d = {\rm rank}\, B$ and $\tilde{d} : = n - 2d$, so that in
the three-dimensional case we have $d=1$ and $\tilde{d}=1$. Note that the
spectrum of the unperturbed Schr\"odinger operator $H_0(B)$ with
magnetic field $B$ is pure point if $\tilde{d}=0$, and is purely
absolutely continuous if  $\tilde{d}\geq 1$. For $\tilde{d}\geq 1$ the unperturbed
operator can be written in appropriate Cartesian coordinates
$(x,y,z) \in {\mathbb R}^n$ with $x,y \in {\mathbb R}^d$ and $z
\in {\mathbb R}^{\tilde{d}}$, as
$$\sum_{j=1}^d \Big( \Big( D_{x_j} + \frac{b_j}{2}y_j \Big)^2+ \Big(
D_{y_j} - \frac{b_j}{2}x_j \Big)^2 \Big) + \sum_{l=1}^{\tilde{d}} D_{z_l}^2,$$
where $B = \sum_{j=1}^d b_j dx_j \wedge dy_j$. We believe that the presence of infinitely many resonances is
typical for the cases $\tilde{d}=1$ and $\tilde{d}=2$. However, the Riemann surfaces where the
resonances are defined, and the eigenvalue counting
functions for the corresponding Toeplitz operators (see \cite{Ra90}) which occur in the
estimates of the resonances, should be of different type in these two
cases. On the other hand, if $\tilde{d} \geq 3$, we expect that the number of the
resonances near any fixed Landau level should be finite. The qualitative pictures in
the cases $\tilde{d}=1$, $\tilde{d}=2$, and $\tilde{d}\geq 3$, should be independent of the rank $2d$
of the non-vanishing magnetic field $B$.
\end{rem}

\section{Perturbation of definite sign}
In this section, we discuss the case $\pm V \geq 0$. We will
obtain an upper bound of the number of resonances near   the
Landau levels outside a semi-axis. Further, for small
perturbations, we prove the existence of a region free of
resonances, and obtain a lower bound on the number of resonances
near  a semi-axis.  In particular, we show that for small positive
perturbations there are no embedded eigenvalues.

In the  definite-sign case, we can summarize  our results by Figure 2.

 Let $V$ have  a definite sign, i.e. let $J = {\rm sign}\,V$
be constant, $J=\pm 1$ when $\pm V\geq 0$. In this case, according
to Proposition \ref{propB.2}, we have
$${\mathcal T}_V(z_{q} (k)) =  \frac{i}{k} J B_q + A(k),$$
with $B_q$  a positive self-adjoint operator independent of $k$,
and $A(k)$  holomorphic near $k=0$ with values in $S_2$. For $iJk
\notin$ sp$(B_q)$, the operator $I + \frac{i}{k} J B_q $ is
invertible with
$$ \|( I + \frac{i}{k} J B_q )^{-1} \| \leq  \frac{|k|}{\sqrt{(J \Im k)^2_+ + |\Re k|^2}},$$
and for $ -\delta J \Im k <  |\Re k|$, the estimate  $ \|( I +
\frac{i}{k} J B_q )^{-1} \| \leq \sqrt{1 + \delta^{-2}}$ holds
uniformly with respect to $k$,  $| k | <s_0$, $ -\delta J \Im k <
|\Re k|$.

 We have
$$ I+{\mathcal T}_V(z_{q} (k)) = \left( I + K(k) \right)\Big(I + \frac{i}{k} J B_q \Big),$$
with
$$K(k): = A(k)\Big(I + \frac{i}{k} J B_q \Big)^{-1}.$$
Note that $K(k) \in S_2$,  and its Hilbert-Schmidt norm is uniformly
bounded with respect to $k$, for  $| k | <s_0$,  $ -\delta J \Im
k <  |\Re k|$. Therefore,
\begin{equation}\label{eq4.1}
{\det}_2\Big(I+{\mathcal T}_V(z_{q} (k))\Big)=
\det\Big(I + \frac{i}{k} J B_q \Big){\det}_2\Big(I+K(k)\Big)
e^{-\Tr({\mathcal T}_V(z_{q} (k))- K(k))} .
\end{equation}
This relation is obtained  by approximating the Hilbert-Schmidt
operator $K$ by a finite-rank operator,  and using the fact that
for a trace-class operator $B$, we have
$\det_2(I+B)=\det(I+B)e^{-\tr\;B}$. We exploit moreover, the fact
that since $B_q$ is a trace class operator (see Corollary
\ref{lemBq}), then  such is $({\mathcal T}_V( z_{q} (k) )-
K(k))=(I+K(k))\frac{i}{k} J B_q$.
\begin{figure}
\begin{center}

\begin{picture}(0,0)%
\includegraphics{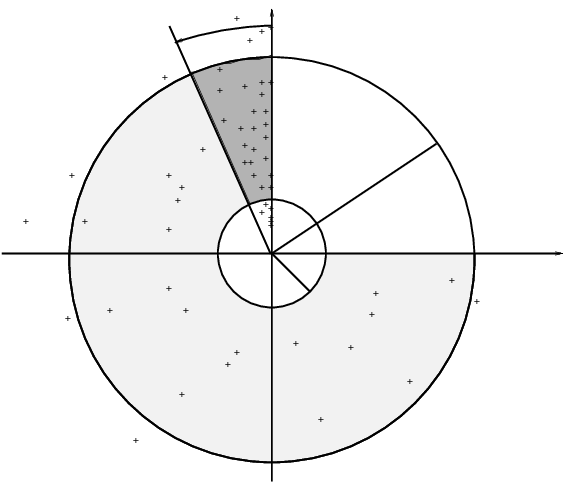}%
\end{picture}%
\setlength{\unitlength}{1066sp}%
\begingroup\makeatletter\ifx\SetFigFont\undefined%
\gdef\SetFigFont#1#2#3#4#5{%
  \reset@font\fontsize{#1}{#2pt}%
  \fontfamily{#3}\fontseries{#4}\fontshape{#5}%
  \selectfont}%
\fi\endgroup%
\begin{picture}(10041,9676)(1168,-9035)
\put(8701,-736){\makebox(0,0)[lb]{\smash{{\SetFigFont{14}{16.8}{\rmdefault}{\mddefault}{\updefault}$V<0$}}}}
\put(6376,-8986){\makebox(0,0)[lb]{\smash{{\SetFigFont{5}{6.0}{\rmdefault}{\mddefault}{\updefault}$ $}}}}
\put(6976,-2536){\makebox(0,0)[lb]{\smash{{\SetFigFont{12}{13.2}{\rmdefault}{\mddefault}{\updefault}$s_0$}}}}
\put(6376,-4336){\makebox(0,0)[lb]{\smash{{\SetFigFont{12}{13.2}{\rmdefault}{\mddefault}{\updefault}$r$}}}}
\put(3976,-5611){\makebox(0,0)[lb]{\smash{{\SetFigFont{14}{16.8}{\rmdefault}{\mddefault}{\updefault}$\bf{C_\theta}$}}}}
\put(5026,-1261){\makebox(0,0)[lb]{\smash{{\SetFigFont{14}{16.8}{\rmdefault}{\mddefault}{\updefault}$\bf{S_\theta}$}}}}
\put(6151,164){\makebox(0,0)[lb]{\smash{{\SetFigFont{12}{13.2}{\rmdefault}{\mddefault}{\updefault}$\hbox{Im} \, k$}}}}
\put(4801,389){\makebox(0,0)[lb]{\smash{{\SetFigFont{12}{13.2}{\rmdefault}{\mddefault}{\updefault}$\theta$}}}}
\put(9826,-3586){\makebox(0,0)[lb]{\smash{{\SetFigFont{12}{13.2}{\rmdefault}{\mddefault}{\updefault}$\hbox{Re} \, k$}}}}
\end{picture}%
\begin{picture}(0,0)%
\includegraphics{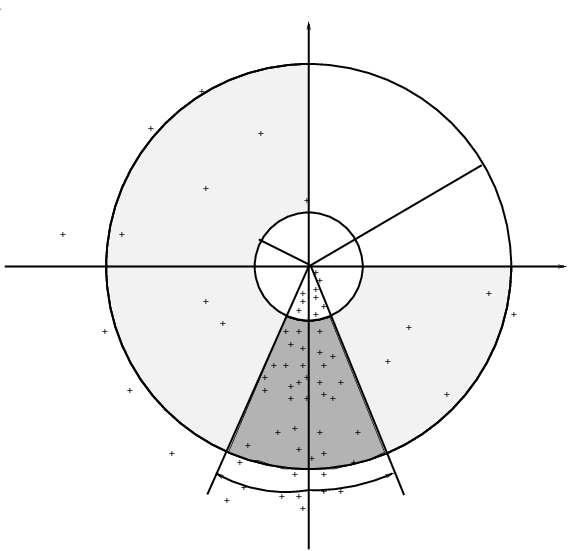}%
\end{picture}%
\setlength{\unitlength}{1066sp}%
\begingroup\makeatletter\ifx\SetFigFont\undefined%
\gdef\SetFigFont#1#2#3#4#5{%
  \reset@font\fontsize{#1}{#2pt}%
  \fontfamily{#3}\fontseries{#4}\fontshape{#5}%
  \selectfont}%
\fi\endgroup%
\begin{picture}(10095,9771)(514,-9019)
\put(7126,-2686){\makebox(0,0)[lb]{\smash{{\SetFigFont{12}{13.2}{\rmdefault}{\mddefault}{\updefault}$s_0$}}}}
\put(5476,-3586){\makebox(0,0)[lb]{\smash{{\SetFigFont{12}{13.2}{\rmdefault}{\mddefault}{\updefault}$r$}}}}
\put(8926,-511){\makebox(0,0)[lb]{\smash{{\SetFigFont{14}{16.8}{\rmdefault}{\mddefault}{\updefault}$V>0$}}}}
\put(5701,-6511){\makebox(0,0)[lb]{\smash{{\SetFigFont{14}{16.8}{\rmdefault}{\mddefault}{\updefault}$\bf{S_\theta}$}}}}
\put(6076,464){\makebox(0,0)[lb]{\smash{{\SetFigFont{12}{13.2}{\rmdefault}{\mddefault}{\updefault}$\hbox{Im} \, k$}}}}
\put(9826,-3661){\makebox(0,0)[lb]{\smash{{\SetFigFont{12}{13.2}{\rmdefault}{\mddefault}{\updefault}${\hbox {Re}} \, k$}}}}
\put(4876,-8461){\makebox(0,0)[lb]{\smash{{\SetFigFont{12}{13.2}{\rmdefault}{\mddefault}{\updefault}$\theta$}}}}
\put(6526,-8461){\makebox(0,0)[lb]{\smash{{\SetFigFont{12}{13.2}{\rmdefault}{\mddefault}{\updefault}$\theta$}}}}
\put(4051,-2386){\makebox(0,0)[lb]{\smash{{\SetFigFont{14}{16.8}{\rmdefault}{\mddefault}{\updefault}${\bf{C_{\theta}}}$}}}}
\end{picture}%
\caption{ {\bf Resonances near a  Landau level for $V$ of definite
sign.} Resonances $z=z_{q} (k)$ are concentrated near the semi
axis $k =- i (\hbox{sgn} V)  ]0, + \infty[$. On the one hand, for
any $\theta$, the number of resonances in $C_\theta$ is bounded by
${\mathcal O}(|\ln r|)$ for $s_0=s_0(\theta)$ sufficiently small
(Proposition \ref{jensensignconstant}). On the other hand, for any
$0<s_0<\sqrt{2b}$ and any $\theta$, there is no resonance of $H_0+
\varepsilon V$  in $C_\theta$ for $\varepsilon\leq \varepsilon_0
(\theta)$ sufficiently small and for compactly supported $V$  we
have lower bound of the number of resonances in $S_\theta$ (see
Theorem \ref{smallV}).}
\end{center}
\label{figres}
\end{figure}

According to (\ref{eq4.1}), for $ | k | <s_0$,  $-\delta J \Im k <
|\Re k|$,  the zeros of  $\det_2(I+ {\mathcal T}_V(z_{q} (k)))$  are the zeros of $\det_2(I+K(k))$.
By  the properties of $K(k)$, $\det_2(I+K(k))= {\mathcal O}(e^{C \|K(k)\|_2^2})={\mathcal O}(1)$, uniformly with respect to $k$. On the other hand,  writing
$$(I+K)^{-1} = (I +   \frac{i}{k} J B_q )\, (I+{\mathcal T}_V)^{-1},$$
and  arguing as in the proof of Proposition \ref{propB.3},
 we find that
$$ \|(I+K)^{-1}\| = {\mathcal O}(|s|^{-1}){\mathcal O}(\delta^{-1})$$
for $\Im k^2> \delta >0$, and for $0<s<|k|<s_0$, uniformly with respect to $(k,s)$.
If $(\lambda_j)_j$ denotes the sequence of eigenvalues of $K(k)$, the above estimate implies that for $\Im k^2> \delta >0$, and for $0<s<|k|<s_0$, we have
\begin{equation}\label{eq4.5}
|1 + \lambda_j|^{-1} =  {\mathcal O}(|s|^{-1}){\mathcal O}(\delta^{-1}).
\end{equation}
 Now, we are able to establish  a lower bound of $\det_2(I+K(k))$. We have
$$\left| \Big({\det}_2(I+K(k))\Big)^{-1}\right| = \left|\det \Big((I+K(k))^{-1}e^{K(k)}\Big)\right|\leq \prod_{|\lambda_j|\leq \frac12} \left|\frac{e^{\lambda_j}}{1 + \lambda_j}\right| \times \prod_{|\lambda_j|> \frac12} \frac{e^{|\lambda_j|}}{|1 + \lambda_j|}.$$
The first product is  uniformly
bounded because $K(k)$ is uniformly
bounded in $S_2$ and we estimate the second product by ${\mathcal O}(e^{ C  (| \ln
  \delta|+ | \ln s|)})$ using  the fact that it involves a finite number of
factors bounded by  ${\mathcal O}(|s|^{-1}){\mathcal O}(\delta^{-1})$ (see
(\ref{eq4.5})).  We get
$$ | {\det}_2(I+K(k))| \geq Ce^{ -C  (| \ln \delta|+ | \ln s|)}$$
for $\Im k^2> \delta >0$, and for $0<s<|k|<s_0$.
Consequently, from  the Jensen inequality (Lemma \ref{lemJensen}), in the case $V$ of definite sign,  we  establish
 upper bounds outside  a neighborhood of   $\{ z_{j} (k) ; \ k \in(-i J) [0, + \infty[\}$:

\begin{prop}\label{jensensignconstant} {\bf Upper bound: special case.}
Suppose that $V$ satisfying (\ref{supexp})
with $\mper >2$, is of definite sign $J$.
For any $\delta >0$,
there exists $s_0>0$, such that for any  $0<r<s_0$ we have
$$\# \{ z = z_{q} (k)  \in \res \cap D_{q}^{*}; \ r< \vert k \vert <
2r  , \  -\delta J \Im k <  |\Re k| \} = {\mathcal O}(| \ln r|).$$
\end{prop}

 In  what follows,  we  prove also that for small perturbations of
definite sign  the
resonances are near $z_{q} (k)$ with $k$ eigenvalues of $-iJB_q$.
In particular, we have a infinite number of resonances  close to
the Landau levels.

 In order to obtain our   lower bound of the counting function of
 resonances, we need the following
 result deduced from Lemma \ref{lem1RW}.

\begin{lem} \label{corRW}
 Let $0\leq W \in L^{\infty} (\R^{2})$ such that
\begin{equation}  \label{sdf}
\ln W (X_{\perp}) \leq - C \langle X_{\perp} \rangle^{2} ,
\end{equation}
for some $C >0$. Let $(\lambda_j)_j$ be the
non-increasing sequence of the non-vanishing eigenvalues of $ p_q W
p_q$, counted with their multiplicity. Then there  exists $\nu>0$ such
that
\begin{equation}
\# \{ j; \; \lambda_j - \lambda_{j+1} > \nu \lambda_j \} = \infty.
\end{equation}
\end{lem}

\begin{proof}
By assumption, one can find a function $U$ which satisfies the
hypotheses of Lemma \ref{lem1RW} with $\beta =1$ such that $W \leq
U$. Then $p_{q} W p_{q} \leq p_{q} U p_{q}$ and
\begin{equation}  \label{dhd}
n_{+} (s; p_{q} W p_{q}) \leq n_ {+} (s; p_{q} U p_{q}) = {\mathcal O}
(\vert \ln s \vert ),
\end{equation}
by the min-max principle.

Let us assume that the set $\{ j; \; \lambda_j - \lambda_{j+1} > \nu \lambda_j \}$ is finite for any $\nu > 0$. Then there exists $j_\nu$ such that for any $j \geq j_\nu$, $\lambda_j - \lambda_{j+1} \leq \nu \lambda_j .$ This implies that for any $j > j_\nu$, $\lambda_j \geq (1-\nu)^{j-j_\nu} \lambda_{j_\nu}.$ In this case for $s$ sufficiently small we would have
$$n_+(s; p_q W p_q) = \# \{j; \, \lambda_j >s \} \geq  \# \{j; \, (1-\nu)^{j-j_\nu} \lambda_{j_\nu}>s \} ,$$
that is
$$n_+(s; p_q W p_q) \geq | \ln(1-\nu)|^{-1} \, |\ln s| - {\mathcal O}_\nu(1).$$
 If we choose $\nu>0$  small enough, this lower bound is in contradiction
with the estimate (\ref{dhd}).
\end{proof}


\begin{thm}\label{smallV}{\bf  Sector free of  resonances,    upper and lower bound.}
 Let $0<s_0<\sqrt{2b}$ and $q \in \N$. Assume $V$ satisfies \eqref{supexp} with $\mper >2$ and is of definite sign $J$.
Then  for any $\delta >0$ there exists $\varepsilon_0 >0$ such
that:

{\rm (i)} for any $\varepsilon \leq \varepsilon_0$,
$H_\varepsilon:=H_0 + \varepsilon V$ has no resonances in $\{z = z_{q}
(k) \in D^{*}_{q} ; \ 0< |k |<s_0, \; -J \Im k \leq \frac{1}{\delta} |\Re
k|\}$.

{\rm (ii)}  there exists $r_0>0$, such that for any
 $0<r<r_0$ and $\varepsilon \leq \varepsilon_{0}$, we have
\begin{equation}  \label{tgb}
\# \{ z = z_{q} (k) \in {\rm Res}(H_\varepsilon) \cap D_{q}^{*} ; \ r < \vert k \vert < 2r \} =
{\mathcal O}( n_+(r, \varepsilon p_q W p_q) - n_+(8r, \varepsilon p_q W p_q)
).
\end{equation}

{\rm (iii)}  if $W$   defined by (\ref{defW})
satisfies (\ref{sdf}),   then  for any $\varepsilon \leq \varepsilon_0$, $H_\varepsilon$ has  an infinite number of
resonances in $\{ z = z_{q} (k) \in D_{q}^{*} ; \ 0< |k |<s_0, \; -J \Im k >
\frac{1}{\delta} |\Re k|\}$.

  More precisely,  there exists a
decreasing sequence $(r_l)_{l\in \N}$ of positive numbers, $r_l
\searrow 0$ such that,
\begin{equation}   \label{Nbinfini}
\begin{aligned}
\# \{ z = z_{q} (k) \in {\rm Res}(H_\varepsilon) \cap D_{q}^{*} ; \ r _{l+1}< \vert
k \vert < r_{l} , \ -J \Im & k > \frac{1}{\delta} |\Re k| \}  \\
&\geq \rank  {\bf 1}_{[ 2r_{l+1} , 2r_{l} ]} ( \varepsilon p_q W p_q ) .
\end{aligned}
\end{equation}

\end{thm}



\begin{proof}
 (i) We have
\begin{equation} \label{rtt}
I+{\mathcal T}_{\varepsilon V} = I+\varepsilon {\mathcal T}_V = I + \frac{i}{k}
\varepsilon J B_q + \varepsilon A(k).
\end{equation}
Since $B_q$ is
self-adjoint and positive, the operator $I + \frac{i}{k} \varepsilon J B_q$ is
invertible for $-J  \Im k < \frac{1}{\delta}
|\Re k|$,  and we have
$$
\|(I + \frac{i}{k} \varepsilon J B_q)^{-1} \| \leq \sqrt{1 +
\delta^{-2}}.
$$
Moreover, for $|k| \leq s_0$ there exists $C>0$
such that $\|A(k)\| \leq C$. Consequently, for $\varepsilon < (C
\sqrt{1 + \delta^{-2}})^{-1}$ and $-J  \Im k \leq
\frac{1}{\delta} |\Re k|$, the operator $I + \frac{i}{k}
\varepsilon J B_q +\varepsilon A(k)$ is invertible and   $z_q(k)$
is not a resonance of $H_0+ \varepsilon V$.

 (ii) We prove this point like Theorem \ref{thmUB}. Let
\begin{align*}
B^{+} := \varepsilon B_q {\bf 1}_{]r/2,4r[}(\varepsilon B_q ) \quad \text{and} \quad
B^{-} := \varepsilon B_q {\bf 1}_{[0 , r/2] \cup [4r , +\infty [}(\varepsilon B_q ).
\end{align*}
For $\frac{2}{3} r < \vert k \vert < \frac{3}{2} r$, the
spectrum of the self-adjoint operator $\frac{1}{\vert k \vert} B^{-}$
is a subset of $[0, \frac{3}{4}] \cup [\frac{8}{3}, + \infty[$. Then
$I + \frac{i}{k} J B^{-}$ is invertible and
\begin{align*}
\Big\Vert \Big( I + \frac{i}{k} J B^{-} \Big)^{-1} \Big\Vert \leq 4 .
\end{align*}
So, if $\varepsilon_{0}$ is small enough and $0 \leq \varepsilon
\leq \varepsilon_{0}$, the operator $I + \frac{i}{k} J B^{-} +
\varepsilon A(k)$ is invertible for $\frac{2}{3} r < \vert k \vert
< \frac{3}{2} r$ with a uniformly bounded inverse.


Using (\ref{rtt}), we can write
\begin{align}
I+{\mathcal T}_{\varepsilon V} =& I +
\frac{i}{k} J B^{+} + \frac{i}{k} J B^{-}+ \varepsilon A(k)  \nonumber  \\
=& \Big( I + \frac{i}{k} J B^{-} + \varepsilon A(k) \Big) \Big( I + \Big(I + \frac{i}{k}
J B^{-} + \varepsilon A(k) \Big)^{-1} \frac{i}{k} J B^{+} \Big),
\end{align}
and, from Proposition \ref{propA.2}, the resonances of $H_{\varepsilon}$, with $\frac{2}{3} r <
\vert k \vert < \frac{3}{2} r$, are the zeros of
\begin{equation}
\widetilde{D} (k,r) := \det \Big( I + \Big(I + \frac{i}{k} J B^{-} + \varepsilon A(k)
\Big)^{-1} \frac{i}{k} J B^{+} \Big).
\end{equation}
Moreover, the multiplicity of the resonance, $\mult (k)$, is  equal to the order of the
zero of $\widetilde{D} (k,r)$. Since $B^{+}/k$ is uniformly bounded, there exists
$C>0$ such that
\begin{equation}  \label{oop}
\vert \widetilde{D} (k,r) \vert \leq \exp \big( C \rank \, {\bf 1}_{] r/2 , 4r [}(
\varepsilon B_q) \big) .
\end{equation}

On the other hand, we have
\begin{align}
\Big(I + \frac{i}{k} J B^{-} + \varepsilon A(k) \Big)^{-1} \frac{i}{k} J
B^{+} =& \Big(I + \Big(I + \frac{i}{k} J B^{-} \Big)^{-1} \varepsilon A(k) \Big)^{-1}
\frac{i}{k} J B^{+}  \nonumber \\
=& \frac{i}{k} J B^{+} + {\mathcal O} (\varepsilon).
\end{align}
For $u \in L^{2} (\R^{3})$, $\frac{2}{3} r < \vert k \vert <
\frac{3}{2} r$ and $k \in \R$, we get
\begin{align}
\Re \Big\langle \Big( I + \Big(I + \frac{i}{k} J B^{-} + \varepsilon A(k)
\Big)^{-1} \frac{i}{k} J B^{+} \Big) u , u \Big\rangle =& \Re
\Big\langle \Big( I + \frac{i}{k} J B^{+} + {\mathcal O} (\varepsilon)
\Big) u , u \Big\rangle  \nonumber \\
=& \Re \big\langle \big( I + {\mathcal O}
(\varepsilon) \big) u , u \big\rangle  \nonumber  \\
\geq& \Vert u \Vert^{2} /2 ,
\end{align}
for $\varepsilon >0$ small enough.  Since we can obtain the same estimate for the adjoint, the
operator $I + (I + \frac{i}{k} J B^{-} + \varepsilon A(k) )^{-1}
\frac{i}{k} J B^{+}$ is invertible for $\frac{2}{3} r < \vert k \vert
< \frac{3}{2} r$, $k \in \R$ with a uniformly bounded
inverse. Then, for such $k$, we get
\begin{align*}
( \widetilde{D} (k,r) )^{-1} =& \det \Big( I - \Big(I + \frac{i}{k} J B^{-} + \varepsilon A(k)
\Big)^{-1} \frac{i}{k} J B^{+} \Big( I + \Big(I + \frac{i}{k} J B^{-} + \varepsilon A(k)
\Big)^{-1} \frac{i}{k} J B^{+} \Big)^{-1} \Big) \\
\leq& \exp \big( C \rank \, {\bf 1}_{] r/2 , 4r [}(
\varepsilon B_q) \big).
\end{align*}
Combining this estimate with (\ref{oop}), the Jensen inequality and
(\ref{nBnW}), we get (\ref{tgb}).

 (iii) According to Lemma \ref{corRW}, there exists
$\nu>0$ and a decreasing sequence $(r_l)_{l\in \N}$ of positive
numbers, $r_l \searrow 0$ such that for any $l \in \N$ we have
$$\hbox{dist}\Big(r_l, \hbox{sp}(B_q)\Big) \geq \nu r_l/2.$$
 Then for any $l \in \N$, there exists a path  (see Figure 3)
$$\widetilde{\Gamma}_l \subset \{ \widetilde{k} \in \C^*; \; |\widetilde{k}| \leq s_0, \; |  \Im \widetilde{k} | \geq
\delta \Re \widetilde{k} , \; r_l \geq \Re \widetilde{k} \geq r_{l+1} \}$$ enclosing the
eigenvalues of $B_q$ contained in the interval $[r_{l+1}, r_l]$, and
such that for $\widetilde{k} \in \widetilde{\Gamma}_l$, the operator $(\widetilde{k}-B_q)$ is invertible with
$$
\|(\widetilde{k}-B_q)^{-1}\|= \sup_{\lambda_j \in \hbox{sp}(B_q)} \frac{1}{|\widetilde{k}-\lambda_j|}
\leq C/|\widetilde{k}|
$$
for some $C=C(\delta, \nu)$, uniformly with respect to $\widetilde{k} \in \widetilde{\Gamma}_l$.

Now, let us consider the path $\Gamma_l:=i\varepsilon J
\widetilde{\Gamma_l}$, and   estimate from below the number of the zeros
of $\det_2( I + \frac{i}{k} \varepsilon J B_q + \varepsilon A(k))$
counted with their multiplicity, enclosed in $\{ z=z_{q} (k) \in
D_{q}^{*}; \ k
\in \Gamma_l\}$.
\begin{figure}
\begin{center}
\begin{picture}(0,0)%
\includegraphics{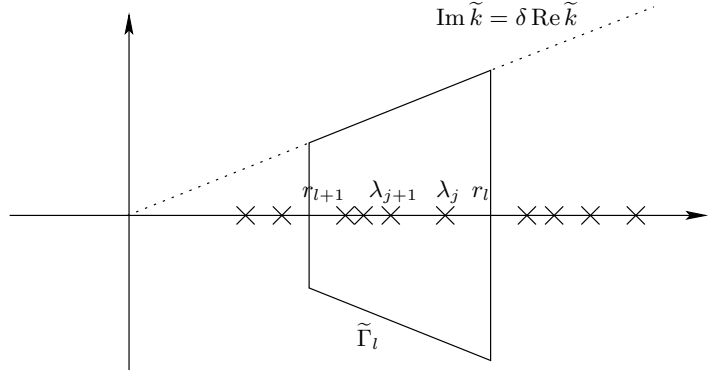}%
\end{picture}%
\setlength{\unitlength}{3947sp}%
\begingroup\makeatletter\ifx\SetFigFont\undefined%
\gdef\SetFigFont#1#2#3#4#5{%
  \reset@font\fontsize{#1}{#2pt}%
  \fontfamily{#3}\fontseries{#4}\fontshape{#5}%
  \selectfont}%
\fi\endgroup%
\begin{picture}(4413,2307)(1189,-3973)
\put(4222,-2875){\makebox(0,0)[rb]{\smash{{\SetFigFont{9}{10.8}{\rmdefault}{\mddefault}{\updefault}$r_{l}$}}}}
\put(3310,-2875){\makebox(0,0)[rb]{\smash{{\SetFigFont{9}{10.8}{\rmdefault}{\mddefault}{\updefault}$r_{l+1}$}}}}
\put(3880,-2875){\makebox(0,0)[lb]{\smash{{\SetFigFont{9}{10.8}{\rmdefault}{\mddefault}{\updefault}$\lambda_{j}$}}}}
\put(3880,-1792){\makebox(0,0)[lb]{\smash{{\SetFigFont{9}{10.8}{\rmdefault}{\mddefault}{\updefault}${\rm Im} \, \widetilde{k} = \delta \, {\rm Re} \, \widetilde{k}$}}}}
\put(3766,-2875){\makebox(0,0)[rb]{\smash{{\SetFigFont{9}{10.8}{\rmdefault}{\mddefault}{\updefault}$\lambda_{j+1}$}}}}
\put(3526,-3811){\makebox(0,0)[rb]{\smash{{\SetFigFont{9}{10.8}{\rmdefault}{\mddefault}{\updefault}$\widetilde{\Gamma}_l$}}}}
\end{picture}%
\caption{The path $\widetilde{\Gamma}_l$}

\end{center}
\label{figgammal}
\end{figure}
 By construction of $\widetilde{\Gamma_l}$, for $k \in \Gamma_l$, the operator $I
+ ik^{-1}\varepsilon JB_q$ is invertible with $\|(I +
ik^{-1}\varepsilon JB_q)^{-1}\|$ $ \leq C(\delta, \nu)$ uniformly
with respect to $k \in \Gamma_l$. Then choosing
$\varepsilon_0$  so small  that
$$\| \varepsilon_0  A(k) (I + ik^{-1}\varepsilon JB_q)^{-1}\|_2 < 1/2,$$
and using that $\det_2(I+A) \leq e^{\|A\|_2/2}$, we obtain  that
for $k \in \Gamma_l$,
$$\Big \vert {\det}_2\Big( I+ \varepsilon_0 A(k)
(I + ik^{-1}\varepsilon JB_q)^{-1} \Big) -1 \Big \vert<1.$$ Applying  the
Rouch\'e theorem we deduce that the number of zeros of  $\det_2( I
+ \frac{i}{k} \varepsilon J B_q + \varepsilon A(k))$ enclosed in $\{
z=z_{q} (k) \in D_{q}^{*} ; \ k \in \Gamma_l\}$ is equal to the number of
zeros of $\det_2( I + \frac{i}{k} \varepsilon J B_q )$ and using
(\ref{nBnW}) it is given by $ n_+(2r_{l+1};  p_q W p_q) -
n_+(2r_{l};  p_q W p_q)$. Since each zero of $\det_2( I +
\frac{i}{k} \varepsilon J B_q + \varepsilon A(k))$ is a resonance, with the same multiplicity, we deduce \eqref{Nbinfini},  and since the sequence $(r_l)_l$ is
infinite, we conclude that the  number of the resonances is
infinite.
\end{proof}

 Since   the embedded eigenvalues in $\R\setminus 2b \N$ are
 resonances $z_{q} (k)$ with $k \in e^{i\{0, \frac{\pi}{2}\}} ]0,\sqrt{2b}[$, a simple consequence of the previous theorem is the absence of embedded eigenvalues in   $]2bq-s_0^2, 2bq[ \cup ]2bq, 2bq+s_0^2[$ for small positive
 $V$ and in $]2bq, 2bq+s_0^2[$ for small negative
 $V$. In fact, by more precise estimates with respect to $q$, for small positive
 $V$, we prove absence of embedded eigenvalues in $\R^+\setminus 2b \N$ and for small negative
 $V$, we obtain information about the localization of the embedded eigenvalues on the left of the Landau levels.
\begin{prop} {\bf Absence of embedded eigenvalues.} \label{cor2}
Assume that $V$ satisfies \eqref{eq1.1} with $\mper >0$ and
$m_3>2$. For a positive potential $V$,   there exists
$\varepsilon_0 >0$ such that for any $\varepsilon \leq
\varepsilon_0$, $H_\varepsilon:=H_0 + \varepsilon V$ has no
embedded eigenvalues  in $\R^+\setminus 2b \N$. For a negative
potential $V$,   there exists $\varepsilon_0 >0$ and $C>0$ such
that for any $\varepsilon \leq \varepsilon_0$, $H_\varepsilon:=H_0
+ \varepsilon V$ has no embedded eigenvalues  in $\R^+\setminus (
2b \N+ ]-\varepsilon C,0[) $.
\end{prop}
\begin{proof} According to  Proposition 2.6 of \cite{BPR}, for $V$ satisfying \eqref{eq1.1}
with $\mper >0$ and $m_3>1$  there exists $C>0$ such that $H_\varepsilon$ has  no embedded
eigenvalues  in $\R^+\setminus ( 2b \N+ ]-\varepsilon C,\varepsilon C[) $.
Then, following the proof of Theorem \ref{smallV} {\rm (i)} (or see proof of Proposition 2.5 of \cite{BPR}),
we have only to check that $\varepsilon_0 >0$ can be chosen independently of $\lambda \in\R^+\setminus 2b \N $
such that for any $\varepsilon \leq \varepsilon_0$ and $\lambda \in\R^+\setminus 2b \N $, $I + \varepsilon
{\mathcal T}_V(\lambda)$ is invertible when $V$ is positive. For negative $V$, we have to choose  $\varepsilon_0 >0$
such that for any $\varepsilon \leq \varepsilon_0$ and $\lambda \in\R^+\setminus ( 2b \N+ ]-b,0[) $, $I + \varepsilon
{\mathcal T}_V(\lambda)$ is invertible.

Let  $\lambda \in\R^+\setminus 2b \N $, then there exists $q \in \N$ and $k \in \C$, $|k|\leq \sqrt{b}$ such that
$\lambda = 2bq + k^2$ ($k\in ]0,\sqrt{b}]$ or $k \in i ]0,\sqrt{b}]$).

We have:
$$
{\mathcal T}_V(\lambda)= J |V|^{\frac12}R_0(\lambda) |V|^{\frac12}= J |V|^{\frac12}\langle x_3\rangle^{\frac{m_3}{2}} \Big( \sum_{j \in \N}  p_j \otimes t_{m_3}(\lambda -2bj) \Big)\langle x_3\rangle^{\frac{m_3}{2}} |V|^{\frac12},$$
where $t_{m_3}$ is the continuous extension of $z \mapsto \langle x_3\rangle^{-\frac{m_3}{2}}(D^2_{x_3}-z)^{-1} \langle x_3\rangle^{-\frac{m_3}{2}}$ from $\Im z > 0$ to $z\in \R\setminus \{0\}$. For $\mu \in  \R\setminus \{0\}$, the integral kernel of $t_{m_3}(\mu)$ is given by:
$$ \langle x_{3} \rangle^{-\frac{m_3}{2}} \;i \frac{ e^{i k |x_3-x_3'|}}{2k}\; \langle x_{3}' \rangle^{-\frac{m_3}{2}},$$
where $k=\sqrt{\mu}$ if $\mu>0$ and $k=i\sqrt{-\mu}$ if $\mu<0$.

It is clear that for $\mu<0$, $\| t_{m_3}(\mu)\| \leq |\mu|^{-1}$
and for $\mu>0$, $\| t_{m_3}(\mu)\| \leq C(m_3) |\mu|^{-\frac12}$
with $C(m_3)=\frac12 \int\langle x_3\rangle^{-m_3} dx_3$ (for more
details, see \cite{BPR}).

From the above estimates, we immediately get:
$$\|  \sum_{j \neq q}  p_j \otimes t_{m_3}(\lambda -2bj)\| \leq \sup_{j \neq q} \| t_{m_3}(\lambda -2bj)\| \leq {\hbox{max}} ( b^{-1},  C(m_3)b^{-\frac12}).$$

Moreover, the series expansion with
respect to $k$ of the kernel of the operator $t_{m_3}$  allows  us to write
$p_q \otimes t_{m_3}(\lambda-2bq)$  as the sum
$$p_q \otimes t_{m_3}(\lambda-2bq)=\frac{i}{k}p_q \otimes \tau  +p_q \otimes \rho(k),$$
where $\tau : L^2(\R) \to L^2(\R)$ is the rank-one operator defined by
$\tau u : = \frac{1}{2}
\big\langle u , \langle . \rangle^{-\frac{m_3}{2}} \big\rangle \langle x_3\rangle^{-\frac{m_3}{2}},$
and $\rho(k)$ is the Hilbert-Schmidt operator with  integral kernel
$$ \langle x_{3} \rangle^{-\frac{m_3}{2}} \;i \frac{ e^{i k |x_3-x_3'|}-1}{2k} \langle x_{3}' \rangle^{-\frac{m_3}{2}}
.$$
Since this integral kernel is bounded by ${\mathcal O}( \langle x_{3} \rangle^{-\frac{m_3}{2}} \; |x_3-x_3'|\;  \langle x_{3}' \rangle^{-\frac{m_3}{2}})$, uniformly with respect to $k$, $|k|\leq b$, it follows that for $m_3>2$, $ p_q \otimes \rho(k)$ is uniformly bounded independently of $q$ and $\lambda$.

 Consequently, for $B:= |V|^{\frac12}\langle x_3\rangle^{\frac{m_3}{2}} \Big(   p_q \otimes \tau \Big)\langle x_3\rangle^{\frac{m_3}{2}} |V|^{\frac12},$ we have
$$\| {\mathcal T}_V(\lambda) -\frac{iJ}{k}B\| \leq M,$$
with $M$ independent of $\lambda$. At last, for $J k \in \R$ or $J
k \in i \R^+$, since $B$ is a positive self-adjoint operator,
$\|(I + i \varepsilon J k^{-1} B )^{-1} \| \leq 1$. Then taking
$\varepsilon_0 < M^{-1}$, for $\varepsilon \leq \varepsilon_0 $,
$I + \varepsilon {\mathcal T}_V(\lambda)$ is invertible for any
$\lambda \in \R^+\setminus  2b \N$ when $J=1$ (i.e. $V \geq 0$)
and for any $ \lambda \in\R^+\setminus ( 2b \N+ ]-b,0[) $ when
$J=-1$ (i.e. $V \leq 0$). This concludes the proof of Proposition
\ref{cor2}.
\end{proof}

\begin{rem}
Further information concerning the localization of
the eigenvalues of the operator $H$ for non-sign-definite
potentials $V$ is contained in \cite[Proposition 2.6]{BPR}.
\end{rem}

\section{Spectral shift function and resonances}

In this section we represent  the derivative of the spectral shift
function (SSF) near  the Landau levels as a sum of  a harmonic
measure related to resonances, and the imaginary part of a
holomorphic function. As in \cite{PZ3}, \cite{BP},  \cite{DZ},
such representation justifies  the Breit-Wigner approximation and
implies  a trace formula.  We deduce also an asymptotic expansion
of the SSF near  a given Landau level; in the case of positive
potentials $V$ which decay slowly enough as $|\xp| \to \infty$
this expansion yields a remainder estimate for the corresponding
asymptotic relations obtained in \cite{FR}.

In the case of a relative trace class perturbation, the SSF is
related to the perturbation determinant by the Krein formula
\begin{equation}\label{krein}
\xi(\lambda)= \frac{ 1}{\pi} \lim_{\varepsilon \rightarrow 0^+} {\rm Arg}\, \det
\Big((H-\lambda -i\varepsilon)(H_0-\lambda -i\varepsilon)^{-1} \Big).
\end{equation}
In our case, $|V|^{\frac12} (H_0+i)^{-1}$ is in the Hilbert-Schmidt class, and
the distribution
\begin{equation}\label{xiprime}
\xi': \; f \in C_0^\infty(\R) \longmapsto  - \, \tr \Big(f(H)-f(H_0)\Big)
\end{equation}
is still well defined, but not the above  perturbation determinant.
Since $V(H_0+i)^{-2}$ is of trace class, we could give a meaning to (\ref{krein}) using meromorphic
extension of the regularized Zeta function (see \cite{Bouclet}), but it will
be more
convenient to introduce the regularized spectral shift function
\begin{equation}\label{defxi2}
\xi_2(\lambda)=  \frac{ 1}{\pi} \lim_{\varepsilon \rightarrow 0^+} {\rm Arg}\, \hbox{det}_2
\Big((H-\lambda -i\varepsilon)(H_0-\lambda -i\varepsilon)^{-1} \Big),
\end{equation}
(see (\ref{eq2.15}) for the definition of $\det_2$)
whose derivative is the following distribution
\begin{equation}\label{xi2prime}
\xi_2': \; f \in C_0^\infty(\R) \longmapsto  - \, \tr \Big(f(H)-f(H_0)-
{\frac{d}{d\varepsilon}f(H_0 + \varepsilon V)}\left|_{\varepsilon=0} \right.\Big)
\end{equation}
(see \cite{kop} or \cite{Bouclet}).  Let us note that in \cite{Bouclet}, these quantities are defined with the opposite sign. We will deduce  the properties of
the SSF from those of the
regularized SSF by using the following    lemma which is
well known for perturbation of the Laplacian (see \cite{kop85},
\cite{BoucletHS}).

\begin{lem}\label{lemlienxixi2}
Let $V$ satisfies \eqref{supexp} with $\mper >2$. On $\R
\setminus 2b \N$, we have
\begin{equation}\label{eqlienxixi2}
\xi' = \xi'_2  + \frac{1}{\pi}\Im  \tr \Big( \partial_z {\mathcal T}_V(.) \Big),
\end{equation}
${\mathcal T}_V(z)$ being defined in Lemma \ref{lemdT}.
\end{lem}

\begin{proof}
According to Lemma \ref{lemdT}, $ \tr ( \partial_z {\mathcal T}_V)$
is analytic on   $\overline{{\mathcal F}_{+}}$.
Then, exploiting (\ref{xiprime}) and (\ref{xi2prime}), we have only to prove:
\begin{equation}\label{topfxixi2}
 \tr \Big({\frac{d}{d\varepsilon}f(H_0 + \varepsilon V)}\left|_{\varepsilon=0} \right. \Big)=
- \frac{1}{\pi} \int_\R f(\lambda) \Im  \tr \Big( \partial_z T(\lambda) \Big) d \lambda,
\end{equation}
for any $f \in C^\infty_0(\R \setminus 2b \N)$.
 By  the Helffer-Sj\"ostrand formula (see for instance \cite{DS}), for  $\widetilde f \in C^\infty_0(\R^2)$ an almost analytic extension of $f$, (i.e.,
$\widetilde f_{\vert \R}=f$ and $\overline\partial_\lambda \widetilde f(\lambda)={\mathcal O}
(|\Im \lambda|^{\infty})$), we have
\begin{equation}\label{formuleHS}
 f(H_0 + \varepsilon V)  =-
\frac{1}{\pi}
 \int_{\mathbb C} \overline\partial {\widetilde f}(z) (z- H_0 -\varepsilon
V)^{-1}  L(dz),
\end{equation}
where $ L(dz)$ denotes the Lebesgue measure on $\C$.
This quantity is  differentiable with respect to $\varepsilon$ and its
derivative at $\varepsilon=0$ is given by
$$
{\frac{d}{d\varepsilon}f(H_0 + \varepsilon V)}\left|_{\varepsilon=0} \right. =
- \frac{1}{\pi}
 \int_{\mathbb C} \overline\partial {\widetilde f}(z) (z- H_0)^{-1}V (z- H_0)^{-1}  L(dz).
$$
Following the proof of Lemma \ref{lemdT}, we check easily that for $\pm \Im z
>0$ the operator $(z- H_0)^{-1}V (z- H_0)^{-1}$ is of trace class with trace
norm bounded by ${\mathcal O} (|\Im z|^{-2})$ and by cyclicity of the trace, for
$\Im z >0$ we have
$$\tr \Big( (z- H_0)^{-1}V (z- H_0)^{-1} \Big) =
\tr \Big( J |V|^{\frac12} (z- H_0)^{-2}  |V|^{\frac12} \Big) =
- \tr  \Big( \partial_z {\mathcal T}_V(z) \Big)$$
and for $\Im z<0$
$$\tr \Big( (z- H_0)^{-1}V (z- H_0)^{-1} \Big)= - \overline{ \tr  \Big( \partial_z {\mathcal T}_V(\overline{z}) \Big)}.$$
 Hence, ${\frac{d}{d\varepsilon}f(H_0 + \varepsilon V)}\left|_{\varepsilon=0} \right.$  is of trace class, and
$$
 \tr \Big({\frac{d}{d\varepsilon}f(H_0 + \varepsilon V)}\left|_{\varepsilon=0} \right. \Big)=
 \frac{1}{\pi} \int_{\Im z >0} \overline\partial {\widetilde f}(z) \tr  \Big(
 \partial_z {\mathcal T}_V(z) \Big)  L(dz) + \frac{1}{\pi} \int_{\Im z <0}
 \overline\partial {\widetilde f}(z)\overline{ \tr  \Big( \partial_z {\mathcal
 T}_V(\overline{z}) \Big)} L(dz).
$$
Then  the Green formula yields  (\ref{topfxixi2}).
\end{proof}

Let $\mathcal W_\pm \subset \subset \Omega_\pm$ be open relatively  compact
subsets of $\pm e^{\pm i ]-2 \theta_0, \varepsilon_0]}]0,2b[$ with
   $\varepsilon_0>0$, $2 \theta_0+ \varepsilon_0<2 \pi$. We assume that these sets
   are  independent of $r$ and that $\mathcal W_\pm$ is simply connected.
Also assume that the intersections between $\Omega_\pm$, $\mathcal
W_\pm$ and $\pm ]0,2b[$ are intervals.  Set  $I_\pm=\mathcal{W}_\pm
\cap \pm$ $]0,2b[$ and choose $0< s_1< \sqrt{\hbox{dist}
(0,\Omega_\pm)}$. In the following, we will identify $2b q + r
\bullet$ with $\pi_{G}^{-1} (2b q + r \bullet ) \cap D_{q}^{*}$, for
$\bullet = \mathcal{W}_\pm$, $\Omega_\pm$, $I_{\pm}$.


\begin{prop} \label{ssf2res}
Let $V$ satisfies (\ref{supexp}) with $\mper >2$. For $\mathcal W_\pm  \subset \subset
\Omega_\pm$ and  $I_\pm$ as above, there exist  functions $g_\pm$
holomorphic in $ \Omega_\pm$,  such
that for $\mu \in 2bq + r I_\pm$,  we have

\begin{align}
 \xi_2'(\mu) = \frac{1}{r\pi} \Im g'_\pm(\frac{\mu-2bq}{r}, r) + \sum_{\fract{w \in \res
\cap 2bq + r \Omega_\pm }{\Im w \neq 0}} \frac{\Im w}{\pi |\mu - w|^2} -
\sum_{w  \in \res \cap 2bq + r I_\pm} \delta (\mu - w) \nonumber \\
- \frac{1}{\pi}\Im  \tr \Big( \partial_z {\mathcal T}_V(\mu) \Big)
\label{eqssf2res}
\end{align}
where  $g_\pm(z, r)$ satisfies the estimate
\begin{equation}
g_\pm(z,r)= {\mathcal O}\left(\left(n_+(s_1\sqrt{r}; p_qWp_q) | \ln r| +
\widetilde{n}_{1}(s_1\sqrt{r}/2) +   \widetilde{n}_{2}(s_1\sqrt{r}/2)\right) \right)
= {\mathcal O}(| \ln r| \; r^{-\frac{1}{\mper}})
\end{equation}
 uniformly with respect to $0<r<r_0$ and
$z \in \mathcal W_\pm$, with $\widetilde n_p$, $p=1$, $2$, defined by
\eqref{ggg4b} or \eqref{deftn}.
\end{prop}

In order to obtain such a representation formula, the first step
is the factorization of the  generalized perturbation determinant.
To this end, we need some complex-analysis results due to J.
Sj\"ostrand, summarized in the following

\begin{prop} \label{propSj} {\rm (\cite{S1}, \cite{S2})}
Let $I$ be a bounded interval in $]0, + \infty[$ and $\Omega$  be a
complex neighborhood of $I$ in  $ e^{i ]-2 \theta_0,
\varepsilon_0]}]0,+ \infty[$ with $\varepsilon_0>0$, $2 \theta_0+
\varepsilon_0<2 \pi$ such that $\Omega \cap \{ \Im z >0; \; \Re z >
0\} \neq \emptyset$. Let $z \mapsto F(z,h)$, $0<h<h_0$,   be a family
of holomorphic functions in $\Omega$, having at  most a finite number
$N(h)$  of zeros in $\Omega$.  Assume  that
$$F(z,h)= {\mathcal O}(1)e^{{\mathcal O}(1)N(h)}, \quad z \in \Omega,$$
  and that for any $\delta >0$ small enough, there exists $C>0$ such that for
$z \in \Omega_{\delta}:=\Omega \cap \{ \Im z >\delta; \; \Re z > 0\ \}$
  we have
  $$| F(z,h) | \geq e^{-C N(h)}.$$
  Then for any simply connected $\widetilde{\Omega} \subset \subset
\Omega$ there exists $g(.,h)$ holomorphic on $\widetilde{\Omega}$ such that
  $$F(z,h) = \prod_{j=1}^{N(h)}(z-z_j)e^{g(z,h)}, \quad \frac{d}{dz}g(z,h) =
      {\mathcal O}(N(h)),
\quad z \in \widetilde{\Omega}.$$
\end{prop}

\begin{proof}[Proof of Proposition \ref{ssf2res}]
Fix $\mathcal W_\pm \subset \subset \Omega_\pm $ ,  and  consider the functions
$$F_\pm : \; z\in \Omega_\pm  \longmapsto D (\sqrt{r }\sqrt{z},\sqrt{r}s_1),$$
   where
\begin{equation}\label{defsqrt}
\sqrt{z} = \left\{\begin{array}{lcc} \sqrt{\rho} e^{i\theta/2} & \hbox{ for} & z=\rho e^{i\theta} \in \Omega_+\\
                                     \sqrt{z} =i \sqrt{\rho} e^{-i\theta/2} & \hbox{ for} &  z=-\rho e^{-i\theta}\in \Omega_-, \end{array}\right.
\end{equation}
 and  $D(k,s)$ is defined by  (\ref{defDk}). The functions $F_\pm$ are holomorphic in $\Omega_\pm$ and ${\widetilde {w}}\in \Omega_\pm$  is a zero of $F_\pm$ if and only if $w= 2bq + {\widetilde {w}}r $ is a resonance of $H$. Then applying Proposition \ref{propSj} to $F=F_+$ and to $F(z)=\overline{F_-(-\overline{z})}$ with $h=r$, $N(r)=n_+(s_1\sqrt{r},p_q W p_q) | \ln r|$, we obtain existence of functions $g_{0,\pm}$  holomorphic in  $\Omega_\pm$
such  that for $z \in \Omega_\pm$,   we have the following factorization:
\begin{equation}
D_\pm (\sqrt{r }\sqrt{z},\sqrt{r}s_1)= \prod_{w \in \res \cap 2bq + r \Omega_\pm}\Big(\frac{z r + 2bq-w}{r} \Big) \; e^{g_{0,\pm}(z ,r)},
\end{equation}
with
\begin{equation}
 \frac{d}{d z} g_{0,\pm}(z,r)= {\mathcal O}(n_+(s_1\sqrt{r},p_q W p_q)| \ln r|),
\end{equation}
 uniformly with respect to $z \in \mathcal W_\pm$.

On the other hand, with the notations of Section 3 (see
Proposition \ref{propB.2}  and the proof of Proposition
\ref{propB.3}), for  $z= z_{q} (\sqrt{r}k)$, $0<s_1<|k|<s_0$ we have
$$ \hbox{det}_2 \Big((H-z)(H_0 -z)^{-1} \Big) = \hbox{det}_2 \Big(I +
{\mathcal T}_V(z) \Big)$$
$$ = D(\sqrt{r }k,\sqrt{r}s_1)) \; \det\Big( \Big(I + \frac{iJ}{\sqrt{r } k}B_q {\bf 1}_{[0,s_1\sqrt{r }/2]}(B_q)+ \widetilde A(k \sqrt{r })\Big) e^{-{\mathcal T}_V(z)}\Big).$$

By  the properties of $ \widetilde A(k)$ (see  the proof of  Proposition
\ref{propB.3}), for $\widetilde K(k) = \frac{iJ}{\sqrt{r } k}B_q {\bf
  1}_{[0,s_1\sqrt{r }/2]}(B_q)+ \widetilde A(k \sqrt{r })$, the difference ${\mathcal T}_V(z)-
\widetilde K(k)$ is a finite-rank operator and as in the proof of (\ref{eq4.1}),
we have
$$
\det\Big( (I + \frac{iJ}{\sqrt{r } k}B_q {\bf 1}_{[0,s_1\sqrt{r }/2]}(B_q)+
\widetilde A(k \sqrt{r }))
e^{-{\mathcal T}_V(z)}\Big)={\det}_2(I+ \widetilde K(k))\; e^{-\Tr({\mathcal T}_V(z)-\widetilde K(k))},$$
where  ${\det}_2(I+ \widetilde K(k))$ is a non-vanishing  holomorphic function, for $0<s_1<|k|<s_0$.
 Since $\widetilde A(k)$ is holomorphic in $S_2$ and
$$
\| \frac{B_q}{s}{\bf 1}_{[0,s]}(B_q)\|_2^2 = - \int_0^s \frac{u^2}{s^2}dn_+(u,B_q)= \widetilde{n}_{2}(s),$$
 we have:
$$\|\widetilde K(k)\|_2^2 = {\mathcal O}(\widetilde{n}_{2}(s_1 \sqrt r/2)),$$
which implies that $|\det_2(1 + \widetilde K(k))| ={\mathcal O} (\exp (\widetilde{n}_{2}(s_1 \sqrt r/2)))$. Using moreover
that $\|\widetilde K(k)\| <1$,
we have
also $|\det_2(1 + \widetilde K(k))|^{-1} =  {\mathcal O}(\exp(\widetilde{n}_{2}(s_1 \sqrt r/2)))$. Then there exists $g_1(.,r)$ holomorphic on $\Omega_\pm$ such
that, $\frac{d}{dz} g_1(z,r)= {\mathcal O}(\widetilde{n}_{2}(s_1 \sqrt
r/2))$, on $\mathcal W_\pm$, and
$${\det}_2(1 + \widetilde K(k))=  e^{g_1(z ,r)}.$$
Then by definition of $\xi_2$ (see (\ref{defxi2})),
 for $\mu = z_{q} (k) \in  2bq + r(\Omega_\pm \cap \R)$ we obtain
\begin{align*}
\xi_2'(\mu) =& \frac{1}{\pi r} \Im \partial_\lambda (g_{0,\pm} + g_1)(\frac{\mu - 2bq}{r}, r) - \sum_{\fract{w \in \res \cap 2bq + r \Omega_\pm}{\Im w \neq 0}} \frac{-\Im w}{\pi |\mu - w|^2} -
\sum_{w  \in \res \cap 2bq + r I_\pm} \delta (\mu - w) \\
&+\frac{1}{\pi} \Im \tr \Big( \frac{1}{2k} \partial_k( \frac{iJ}{ k}B_q {\bf 1}_{[0,s_1\sqrt{r }/2]}(B_q)+ \widetilde A(k)) - \partial_z {\mathcal T}_V (\mu + i0)\Big),
\end{align*}
where
$$
k=\left\{
\begin{array} {l}
\sqrt{\mu-2bq} \quad {\rm if} \quad \mu-2bq>0,\\
i\sqrt{2bq-\mu}  \quad {\rm if} \quad \mu-2bq<0.
\end{array}
\right.
$$
According to Lemma \ref{lemdT} and since $B_q \in S_1$, the  operators
$\partial_z{\mathcal T}_V(z)$ and  $\partial_k\widetilde
A(k)=\partial_kA(k)=\partial_k({\mathcal T}_V(z_{q} (k))- \frac{i}{k}JB_q) $ are of trace class. The trace of $\partial_kA(k)$ is given by the integral of its kernel on the diagonal:
\begin{equation}\label{trdA}
\tr \Big( \partial_kA(k)\Big)= \frac{2kb}{8\pi} \int_{\R^3}V(x) dx \Big( \sum_{j>q} (2b(j-q)-k^2)^{-3/2}- i \sum_{j<q} (k^2 + 2b(q-j))^{-3/2}\Big),
\end{equation}
and by definition of $\widetilde n_1$ (see \eqref{deftn})
$$\tr \Big( \frac{1}{2k} \partial_k( \frac{iJ}{ k}B_q
    {\bf 1}_{[0,s_1\sqrt{r }/2]}(B_q)\Big)= -\frac{iJs_1\sqrt{r }}{8k^3}\widetilde{n}_{1}(s_1 \sqrt r/2).$$
Then, we conclude the proof of  Proposition \ref{ssf2res} with  $g_\pm = g_{0,\pm}+g_1 + g_2$  taking
$$g_2(z,r) :=  \frac{ib}{4\pi} \int_{\R^3}V(x) dx  \sum_{j<q} (zr + 2b(q-j))^{-1/2}+ \frac{iJs_1}{4 \sqrt{z}} \widetilde{n}_{1}(s_1 \sqrt r/2),$$
where $\sqrt{z}$ is defined on $\Omega_\pm$ by \eqref{defsqrt}.
\end{proof}

In the definite-sign case ($J= \hbox{sign}\,V$) we can specify the
representation of the regularized spectral shift function when for $z
= z_{q} (k)$,  the operator $1+ \frac{iJ}{k} B_q$ is invertible, that is for ${\rm Arg}\, k \neq -J \pi/2$. Then we consider  $\mathcal W_\pm \subset \subset  \Omega_\pm$  open relatively   compact subsets  of $\pm e^{\pm i ]-2 \theta_0, \varepsilon_0]}]0,2b[$ with $\varepsilon_0>0$, $2 \theta_0+ \varepsilon_0<2 \pi$ as above, and we have the assumption that
\begin{equation}\label{hyptheta}
-J \frac{\pi}{2} \notin (\frac{\pi}{2})_\mp \pm [-\theta_0, \varepsilon_0/2],
\end{equation}
where $(\frac{\pi}{2})_-=0$ and $(\frac{\pi}{2})_+=\frac{\pi}{2}$. The main restriction is in the case "-" for $V\leq 0$ (i.e. $J=-1$),  where we can consider $\Omega_- \subset  -  e^{- i ]-2 \theta_0, \varepsilon_0]}]0,2b[$, only with $\theta_0<0$, that is  where there are no resonances! So in the definite-sign case
 we discuss the three following situations:

For $V \geq 0$, we consider   $\mathcal W_+ \subset \subset  \Omega_+ \subset$$e^{ i ]-2 \theta_0, \varepsilon_0]}]0,2b[$ with $\varepsilon_0>0$, $2 \theta_0< \pi$ or $\mathcal W_- \subset \subset  \Omega_-\subset $$- e^{- i ]-2 \theta_0, \varepsilon_0]}]0,2b[$ with $\varepsilon_0>0$, $2 \theta_0+ \varepsilon_0<2 \pi$.

For $V \leq 0$    we consider  $\mathcal W_+ \subset \subset \Omega_+ \subset$$ e^{ i ]-2 \theta_0, \varepsilon_0]}]0,2b[$ with $\varepsilon_0>0$, $2 \theta_0+ \varepsilon_0<2 \pi$.

\begin{prop} \label{ssf2resbis}
Assume $V$ satisfies (\ref{supexp}) with $\mper >2$, and is of definite sign $J = {\rm sign}V$.
Let $\mathcal W_\pm \subset \subset \Omega_\pm$ open relatively  compact
subsets  of $\pm e^{\pm i ]-2 \theta_0, \varepsilon_0]}]0,2b[$as above and let
         the interval $I_\pm = \mathcal W_\pm \cap \pm$ $]0,2b[$.
Assume $\theta_0$ satisfies (\ref{hyptheta}).

Then for $\lambda=\mu - 2bq \in r I_\pm$, the representation (\ref{eqssf2res}) holds with
$$\frac{1}{r} \Im  g'_\pm(\frac{\lambda}{r},r) = \frac{1}{r}\Im \widetilde g'_\pm(\frac{\lambda}{r},r) +\Im \widetilde  g'_{1,\pm}(\lambda) + {\bf 1}_{(0,2b)}(\lambda)J \, \Phi^\prime (\lambda),$$
where
$$\Phi(\lambda):=\tr \Big( \arctan \frac{B_q}{\sqrt {\lambda}}\Big)=\tr \Big( \arctan \frac{p_qWp_q}{2\sqrt {\lambda}}\Big),$$ $z \mapsto  \widetilde g_\pm(z,r)$ is holomorphic in $ \Omega_\pm$ and satisfies
\begin{equation}
\widetilde{g}_\pm(z, r)= {\mathcal O}(| \ln r|)
\end{equation}
 uniformly with respect to $0<r<r_0$ and
$z \in  \mathcal W_\pm$   while the function $\widetilde g_{1,\pm}$ is holomorphic on $\pm e^{\pm i ]-2 \theta_0, \varepsilon_0]}]0,2b[$ and for $z \in \pm e^{ \pm i ]-2 \theta_0, \varepsilon_0]}]0,2b[$, there exists $C_{\theta_0}$ such that:
$$|\widetilde g_{1,\pm}(z)| \leq C_{\theta_0} \sigma_2(\sqrt {|z|})^{\frac12},$$
$\sigma_2$ being defined in Corollary \ref{lemBq}.


\end{prop}

\begin{proof}
 With   the notations of  Section 3 (see Proposition \ref{propB.2}), and according to relation (\ref{eq4.1}), for  $z= z_q(k)$, $0<s_1 \sqrt{r}<|k|<s_0$, $-J  \Im k < \frac{1}{\delta}
|\Re k|$  we have
$$ \hbox{det}_2 \Big((H-z)(H_0 -z)^{-1} \Big) = \hbox{det}_2 \Big(I +
{\mathcal T}_V(z) \Big)$$
$$= {\det}_2\Big(I+K(k)\Big) \det\Big(I + \frac{i}{k} J B_q \Big) e^{-\Tr({\mathcal T}_V(z)-K(k))}.$$
with $K(k) = A(k)\Big(I + \frac{i}{k} J B_q \Big)^{-1}$.
As in the    proof of Proposition \ref{ssf2res},   applying Proposition
\ref{propSj} and  the results of Section 4, we obtain existence of functions $\widetilde g_\pm$  holomorphic in  $\Omega_\pm$
such  that for $z \in \Omega_\pm$,   we have the following factorization:
\begin{equation}
{\det}_2\Big(I+K(\sqrt{r }\sqrt{z})\Big)= \prod_{w \in \res \cap 2bq + r \Omega_\pm}\Big(\frac{z r + 2bq-w}{r} \Big) \; e^{\widetilde g_\pm (z ,r)},
\end{equation}
with
\begin{equation}
 \frac{d}{d\lambda} \widetilde g_\pm (z,r)= {\mathcal O}(| \ln r|),
\end{equation}
 uniformly with respect to $z \in \mathcal W_\pm$.
Then
by definition of $\xi_2$ (see (\ref{defxi2})),
for $\mu = z_{q} (k) \in  2bq + r I_{\pm}$ we obtain:
\begin{align*}
\xi_2'(\mu) =& \frac{ 1}{\pi r} \Im \widetilde g'_\pm(\frac{\mu - 2bq}{r}, r) + \sum_{\fract{w \in \res
\cap 2bq + r \Omega_\pm}{\Im w \neq 0}} \frac{\Im w}{\pi |\mu - w|^2} - \sum_{w  \in \res \cap 2bq + r I_\pm} \delta (\mu - w)   \\
&+\frac{1}{2k \pi} \Im \tr \Big( (I + \frac{iJ}{ k}B_q)^{-1} \;  \partial_k( \frac{iJ}{ k}B_q ) \Big) - \frac{1}{\pi} \Tr \Big(\partial_z  {\mathcal T}_V (\mu) -\frac{1}{2k}  \partial_k K(k)\Big),
\end{align*}
where
$$
k=\left\{
\begin{array} {l}
\sqrt{\mu-2bq} \quad {\rm if} \quad \mu-2bq>0,\\
i\sqrt{2bq-\mu}  \quad {\rm if} \quad \mu-2bq<0.
\end{array}
\right.
$$

On the other hand, Lemma \ref{lemdT} and Corollary \ref{lemBq} entails that the operators $\partial_z {\mathcal T}_V(\mu)$ and
$$ \partial_k K(k) =  \partial_k A(k)-  \partial_k \Big( A(k) \frac{iJ}{ k}B_q (I + \frac{iJ}{ k}B_q)^{-1}\Big),$$
are trace-class. Moreover, from \eqref{trdA} and since $A$ is holomorphic in $S_2$, we have
$$ \Im \frac{1}{2k} \Tr ( \partial_k  K(k) )= \Im \frac{1}{2k} \partial_k \Big( \widetilde g_{1,\pm} (k^2)\Big),$$
with $ \widetilde g_{1,\pm}$ being the holomorphic function:
$$ \widetilde g_{1,\pm}(z) :=  \frac{ib}{4\pi} \int_{\R^3}V(x) dx  \sum_{j<q} (z + 2b(q-j))^{-1/2}  - \Tr \Big(A(\sqrt{z}) \, \frac{iJ}{ \sqrt{z}} B_q (I + \frac{iJ}{ \sqrt{z} }B_q)^{-1}\Big),$$
which satisfies the claimed estimates thanks to \eqref{eqesti}.
Finally, we have
\begin{align*}
\frac{1}{2k} \Im \, \tr \Big( (I + \frac{iJ}{ k}B_q)^{-1} \partial_k( \frac{iJ}{ k}B_q )\Big) &=
\frac{-1}{2 k^2} \Im \tr \Big(\frac{iJ}{ k}B_q  (I + \frac{iJ}{ k}B_q)^{-1}\Big)  \\
&=\left \{
\begin{aligned}
&0 && \hbox{for } J k \in i \R^+\\
&\frac{1}{2 k^2}  \tr \Big(\frac{J}{ k}B_q  (I + \frac{B_q^2}{ k^2})^{-1}\Big)=
J \Phi^\prime (k^2) \quad && \hbox{for } k \in  \R \, ,
\end{aligned} \right.
\end{align*}
which prove Proposition \ref{ssf2resbis}.
\end{proof}

Combining Lemma \ref{lemlienxixi2}, Proposition \ref{ssf2res} and Proposition \ref{ssf2resbis}, we
deduce Breit-Wigner approximation of the SSF:

\begin{thm}\label{thmBW} {\bf Breit-Wigner approximation.}
Assume $V$ satisfies (\ref{supexp}) with $\mper >2$. Let $\mathcal W_\pm \subset \subset \Omega_\pm$ be open relatively compacts subsets of $\pm e^{\pm i ]-2 \theta_0, \varepsilon_0]}]0,2b[$ as before Proposition \ref{ssf2res}   and let $0< s_1<\sqrt{ \hbox{dist}(\Omega_\pm, 0)}$.   Then there exists $r_0>0$ and
  functions $g_\pm$ holomorphic in  $  \Omega_\pm$,  such
that for $\mu \in 2bq + r I_\pm$,  we have

\begin{equation} \label{ssfres}
 \xi'(\mu) =  \frac{1}{r\pi} \Im g'_\pm(\frac{\mu-2bq}{r}, r) + \sum_{\fract{w \in \res
\cap 2bq + r \Omega_\pm}{\Im w \neq 0}} \frac{\Im w}{\pi |\mu - w|^2} -
\sum_{w  \in \res \cap 2bq + r I_\pm} \delta (\mu - w)
\end{equation}
where  $g_\pm(z, r)$ satisfies the estimate
\begin{equation}\label{estig}
g_\pm(z,r)= {\mathcal O}\left(n_+(s_1\sqrt{r}; p_qWp_q) | \ln r| +
\widetilde{n}_{1}(s_1\sqrt{r}/2) +   \widetilde{n}_{2}(s_1\sqrt{r}/2)\right)
= {\mathcal O}(| \ln r| \; r^{-\frac{1}{\mper}})
\end{equation}
 uniformly with respect to $0<r<r_0$ and
$z \in   \mathcal W_\pm$.

Moreover for potentials of definite sign, $J:={\rm sign}\,V$, assuming that
                   $\theta_0$ satisfies
$$-J \frac{\pi}{2} \notin (\frac{\pi}{2})_\mp \pm [-\theta_0, \varepsilon_0/2],$$
($(\frac{\pi}{2})_-=0, $ $(\frac{\pi}{2})_+=\frac{\pi}{2}$) for $\lambda \in r
                   I_\pm$ we have

$$\frac{1}{r} \Im  g'_\pm(\frac{\lambda}{r},r) = \frac{1}{r}\Im \widetilde g'_\pm(\frac{\lambda}{r},r) +\Im \widetilde  g'_{1,\pm}(\lambda) + {\bf 1}_{(0,2b)}(\lambda)J \,  \Phi^\prime (\lambda),$$
where
$$\Phi(\lambda):= \tr \Big( \arctan \frac{B_q}{\sqrt {\lambda}}\Big)= \tr \Big( \arctan \frac{p_qWp_q}{2\sqrt {\lambda}}\Big),$$ $z \mapsto  \widetilde g_\pm(z,r)$ is holomorphic in $ \Omega_\pm$ and satisfies
\begin{equation}
\widetilde{g}_\pm(z, r)= {\mathcal O}(| \ln r|)
\end{equation}
 uniformly with respect to $0<r<r_0$ and
$z \in  \mathcal W_\pm$. The function $\widetilde g_{1,\pm}$ is holomorphic on
 $\pm e^{\pm i ]-2 \theta_0, \varepsilon_0]}$ $]0,2b[$ and for $z \in \pm e^{ \pm i ]-2 \theta_0, \varepsilon_0]}]0,2b[$, there exists $C_{\theta_0}$ such that:
\begin{equation}\label{estitg1}
|\widetilde g_{1,\pm}(z)| \leq C_{\theta_0} \sigma_2(\sqrt {|z|})^{\frac12},
\end{equation}
$\sigma_2$ being defined in Corollary \ref{lemBq}.
\end{thm}


The following corollary describes the asymptotic behavior of the SSF on the
right of a given Landau level.

\begin{cor}\label{thmFR} {\bf Singularities at the Landau levels.}
Assume  that $V$ satisfies (\ref{supexp}) with $\mper >2$ and is of
definite sign $J={\rm sign}\,V$. Then the asymptotic relation
\begin{equation}\label{asympFR}
\xi(2bq +\lambda)= \frac{J}{\pi} \Phi(\lambda) + {\mathcal O}\left(  \Phi(\lambda)^{\frac12} \right) + {\mathcal O}\left( | \ln \lambda|^{2}  \right),
\end{equation}
holds as $\lambda \searrow 0$.

 \end{cor}

\begin{proof}
Let us apply Theorem \ref{thmBW} on intervals $2bq+ r_n[1,2]$, with  $r_n =
\lambda 2^n$, $\lambda >0$. For $\mu \in 2bq + r_n[1,2]$ and
$\Omega_+$ a complex neighborhood of $[1,2]$, we have
\begin{align}
\xi^\prime (\mu)=& \frac{1}{r_n \pi} \Im \widetilde g_\pm^\prime \left(\frac{\mu-2bq}{r_n}, r_n\right)
 + \sum_{\fract{w \in \res \cap 2bq + r _n\Omega_+}{\Im w \neq 0}} \frac{\Im w}{\pi |\mu - w|^2} -
\sum_{w  \in \res \cap 2bq + r_n[1,2]} \delta (\mu - w)   \nonumber  \\
&+\frac{1}{\pi} \left( J \Phi^\prime + \Im \widetilde  g'_{1,\pm}\right) (\mu - 2bq).  \label{eq5.20}
\end{align}
Using that $\int_\R \frac{-\Im w}{\pi |\mu - w|^2} d \mu  \leq 1$ and that
$2bq + r _n\Omega_+$ contains at the most ${\mathcal O}(| \ln r_n |)$
resonances (see Proposition \ref{jensensignconstant}), integration of
(\ref{eq5.20}) on $2bq + r_n[1,2]$ yields
\begin{equation}\label{eq5.21}
\xi(2bq + r_{n+1}) -\xi(2bq + r_{n}) = \frac{1}{\pi} \left[\Im \widetilde g_\pm(\,.\, , r_n)\right]_{1}^2 + {\mathcal O}(| \ln r_n |)
+\frac{1}{\pi}  \left[J \Phi + \Im \widetilde  g_{1,\pm} \right]_{r_n}^{r_{n+1}}.
\end{equation}
Let $N \in \N$ such that $\frac{b}{2} \leq \lambda 2^{N+1} \leq b$ (then $N =
{\mathcal O}(| \ln \lambda |)$). Since $\xi$, $\Phi$ and $ \widetilde g_{1,\pm}$ are uniformly bounded
on $2bq + b[1/2,1]$ ($b$ the fixed strength of the  magnetic field), and since $\widetilde g_\pm(\,.\, , r_n) = {\mathcal O}(| \ln r_n |)$, taking the sum of (\ref{eq5.21}) from $n=0$ to $n=N$, we have:
$$\xi(2bq + \lambda) = \frac{J}{\pi} \Phi(\lambda) +\frac{1}{\pi} \Im \widetilde g_{1,\pm}(\lambda)+ \sum_{n=0}^N{\mathcal O}(| \ln (2^n \lambda) |) + {\mathcal O}(1).$$
Using \eqref{estitg1} and exploiting that  $N = {\mathcal O}(| \ln \lambda |)$, we obtain existence of $C>0$ such that
$$\left| \xi(2bq + \lambda) - \frac{J}{\pi} \Phi(\lambda) \right| \leq C | \ln \lambda |^2+ C \sigma_2(\sqrt {\lambda})^{\frac12}.
$$
Then   Corollary \ref{thmFR} follows from the   elementary inequality
$$\frac{u^2}{1+ u^2} \leq \arctan u, \quad u \geq 0,$$
which implies    $\sigma_2(\sqrt{\lambda}) \leq \Phi(\lambda)$.
\end{proof}
Let us compare our results with those of \cite{FR} where the singularities of
the SSF at a given Landau level were investigated.

If $W$ satisfies the assumptions of Lemmas \ref{lemRa90}, \ref{lem1RW} or
\ref{lem2RW} (plus some generic technical assumptions in the case of a rapid
decay), then it is shown in \cite{FR} that
\begin{equation} \label{ggg14}
\xi(2bq +\lambda)= \frac{J}{\pi} \Phi(\lambda) \, (1 + o(1)), \quad \lambda
\searrow 0.
\end{equation}
In the case of slowly decaying $W$ satisfying the assumptions of Lemma
\ref{lemRa90}, we see that  \eqref{asympFR} provides a remainder estimate of
asymptotic relation \eqref{ggg14}. In the case of rapidly decaying $W$ (see
Lemma \ref{lem1RW} or
Lemma \ref{lem2RW}) we have
$$
\Phi(\lambda) = \frac{1}{2} \varphi_{\beta}(\lambda) (1 + o(1)), \quad \lambda
\searrow 0.
$$

Hence, in the case $\beta \in ]0, \frac{1}{2}[$ relation \eqref{asympFR} again  provides a remainder estimate of
\eqref{ggg14}. However, in the case $\beta \in [\frac{1}{2},\infty]$ it does not even
recover the first asymptotic term of \eqref{ggg14}. Note also that in
\cite{FR} the decay of $V$ is assumed to be isotropic in all three directions
while here we  assume that $V$ is super-exponentially decaying with
respect to $x_3$. \\
On the left of Landau level, for $V\geq 0$, the results of \cite{FR} imply
$\xi(2bq -\lambda)= {\mathcal O}(1)$. Here our estimates are not accurate to
see this. For $V\leq 0$, it is shown in \cite{FR} that
$$
\xi(2bq - \lambda) = n_+(2\sqrt{\lambda}; p_qWp_q) (1 + o(1)),  \quad \lambda
\searrow 0.
$$
In this case, we have only general representation formula (\ref{ssfres}) with estimate (\ref{estig}).



As in   \cite{PZ3}, (or \cite{BP}), from Theorem \ref{thmBW} we deduce also the following trace formula.

\begin{cor}{\bf Trace formula.}
 Let  $\mathcal W_\pm \subset \subset \Omega_\pm$ be as in Theorem \ref{thmBW}.
Suppose that $f_\pm$ is holomorphic on a neighborhood of  $\Omega_\pm$
and that $\psi_\pm \in C_0^{\infty}( \Omega_\pm\cap \R)$ satisfies $\psi_\pm(\lambda)=1$ near ${\mathcal W}_\pm$.
 Then under the assumptions of
Theorem \ref{thmBW} we have the following trace formula
\begin{equation}
\label{eq:4.1}
 \tr \Big( (\psi_\pm  f _\pm) (\frac{ H - 2bq}{r})- ( \psi_\pm f_\pm)( \frac{ H_0 - 2bq }{r} )\Big)  =
\sum_{ w \in \res \cap 2b q + r {\mathcal W}_\pm} f_\pm( \frac{ w - 2bq }{r})
+ E_{ f _\pm, \psi_\pm}  ( r)
\end{equation}
with
\[ | E_{ f _\pm, \psi_\pm}  ( r)| \leq M( \psi_\pm )
 {\rm{sup} } \;  \{ |f_\pm ( z ) | \; : \; z \in \Omega_\pm \setminus{\mathcal W}_\pm \,, \ \Im z \leq 0 \} \times N_q(r)  \,. \]
where  $N_q(r)=n_+(s_1\sqrt{r}; p_qWp_q) | \ln r| +
  \widetilde{n}_{1}(s_1\sqrt{r}/2) +   \widetilde{n}_{2}(s_1\sqrt{r}/2) = {\mathcal O}(| \ln r| \; r^{-\frac{1}{\mper}})$.
\end{cor}
\begin{proof}
In the proof we omit the subscript $\pm$.
Choose an almost analytic extension $\widetilde{\psi}$ of $\psi$ so that
$\widetilde{\psi} \in  C_0^{\infty} ( \Omega)$, $\widetilde{\psi} = 1$ on
${\mathcal W}$ and
$$ {\rm{supp}}\:\: \overline{\partial}_z \widetilde{\psi}\subset \Omega \setminus {\mathcal W}.
$$
Applying Theorem \ref{thmBW}, we have
\begin{align*}
\tr \Big( (\psi  f ) & (\frac{ H - 2bq}{r})- ( \psi f)( \frac{ H_0 - 2bq }{r} )\Big)  = - \Big\langle \xi'(\lambda), (\psi f)(\frac{ \lambda - 2bq }{r}) \Big\rangle  \\
=& \sum_{w \in \res \cap 2b q + r \supp \psi}
(\psi f)( \frac{ w - 2bq }{r})  - \frac{1}{r \pi }\int (\psi f)(\frac{ \lambda - 2bq }{r})   \Im g'(\frac{\lambda-2bq}{r}, r)
 d\lambda  \\
&+ \frac{1}{2 \pi i} \int (\psi f)(\frac{ \lambda - 2bq }{r})  \sum_{\fract{w \in
\res \cap 2bq + r \supp \widetilde{\psi}}{\Im w \neq 0}}
\Bigl(\frac{1}{\lambda -w } - \frac{1}{\lambda -\overline{w} }
\Bigr)d\lambda \,.
\end{align*}
The integral involving $g'$ can be estimated using
(\ref{estig}) on $ {\rm{supp}}\: \widetilde{\psi}$. For
the integral related to the resonances,
we apply  the Green formula and we get the term
\begin{align*}
\sum_{w \in \res ,  \Im w \neq 0}&(\widetilde{\psi}f)(\frac{ w - 2bq }{r}) \\
+ \frac{1}{r\pi} & \int_{\C_{-}} (\overline{\partial}_z\widetilde{\psi})(\frac{ z - 2bq }{r})
f(\frac{ z - 2bq }{r}) \sum_{\fract{w \in \res \cap 2b q + r \supp \widetilde{\psi}}{\Im w \neq 0}} \Bigl(\frac{1}{z - \overline{w}} - \frac{1}{z - w} \Bigr) L(dz) \,.
\end{align*}
 We apply the inequality
\begin{equation*}
\int_{\Omega_1} \frac{1}{|z - w|}  L(dz) \leq 2 \sqrt{2 \pi
|\Omega_1|} ,
\end{equation*}
and the upper bound of the resonances in $\Omega$
contained in Theorem  \ref{thmUB},
to obtain the result.
\end{proof}

{\bf Acknowledgments.}  The authors would like to thank the referee
for the  
helpful suggestions concerning the descriptions of the subsets
of ${\mathcal M}$. Jean-Fran\c{c}ois Bony and Vincent Bruneau
are partially supported by the French ANR  Grant no. JC0546063.
Vincent Bruneau and Georgi Raikov are partially supported by the Chilean Science
Foundation {\em Fondecyt} under Grants 7050263 and 1050716.

\end{document}